%% file: paper.tex
\def\useplain{1}

\ifx \useplain\undefined
	\documentclass[mnsc]{../style/informs3aa}
	\usepackage{style_set/informs_set}
	\input{title_info/informs_title}

\else
	\documentclass{article}
	\usepackage[margin=1.5in]{geometry}
		\usepackage{style_set/plain_set}
	\input{title_info/plain_title}

\fi

\usepackage{eie_pckg}
\usepackage{eie}

\graphicspath{{./figures/}}

\begin{document}

\ifx \useplain\undefined
\else
\maketitle
\fi

\def\atxt{
Classical approaches to experimental design assume that intervening on one unit does not affect other units. There are many important settings, however, where this non-interference assumption does not hold, as when running experiments on supply-side incentives on a ride-sharing platform or subsidies in an energy marketplace. In this paper, we introduce a new approach to experimental design in large-scale stochastic systems with considerable cross-unit interference, under an assumption that the interference is structured enough that it can be captured via mean-field modeling. Our approach enables us to accurately estimate the effect of small changes to system parameters by combining unobstrusive randomization with lightweight modeling, all while remaining in equilibrium. We can then use these estimates to optimize the system by gradient descent. Concretely, we focus on the problem of a platform that seeks to optimize supply-side payments $p$ in a centralized marketplace where different suppliers interact via their effects on the overall supply-demand equilibrium, and show that our approach enables the platform to optimize $p$ in large systems using vanishingly small perturbations.

\noindent
\emph{Keywords}: experimental design, interference, mean-field model, stochastic system. }

\ifx \useplain\undefined
\ABSTRACT{\atxt}
\else
\begin{abstract}
\atxt
\end{abstract}
\fi


\ifx \useplain\undefined
\maketitle
\fi


\section{Introduction}

\ifx \useplain\undefined
 Randomized controlled trials
\else
Randomized controlled trials\blfootnote{\hspace{-5.3mm}This
work was partially supported by a seed grant from the Stanford Global Climate
and Energy Project and a Facebook Faculty Award.}
\fi
are widely used to guide decision making across different domains, ranging from
classical industrial and agricultural applications \citep{fisher1935design} to developmental economics \citep{banerjee2011poor} and the modern technology sector \citep{athey2019economists,kohavi2009controlled,tang2010overlapping}.
In its most basic form, a randomized trial aims to assess the expected effectiveness of a set of interventions on a population by selecting a small but representative sub-population of units and assigning to each unit a randomly chosen intervention.
For example, in a medical trial, the decision maker may want to compare the effectiveness of a new experimental drug with the current standard of care. To do so they select a set of patients, and randomly assign some fraction to the new treatment while others are given the control condition (i.e., current standard of care). The drug is then assessed by comparing the outcomes of treated and control patients.
Similar randomized experiments  are popular with technology companies, where they are often referred to as A/B tests. In this context, a company would select a small population of its users and expose them to different randomly generated designs; the best design that emerges from the experiment is then deployed to the entire user base at large. 

When interpreting the results from randomized trials, it is common to make a ``no interference'' assumption, whereby
we assume that the intervention assigned to any given unit does not affect observed outcomes for other units \citep{imbens2015causal};
for example, in our medical example, we might assume that giving the experimental treatment to some patients does not
affect outcomes for the control patients who are still receiving standard care.
Such a lack of interference plays a key role in enabling us to use randomized trials to understand the effect of large scale policy interventions,
as it implies that any effects observed by experimenting on a representative sub-population should also hold when the same interventions
are applied to the overall population at large.
However, this non-interference  assumption is violated in many important applications,
and randomized trials can lead to highly misleading conclusions in the presence of cross-unit interference.
We illustrate this problem below using an example of \citet{heckman1998general}.

\begin{exam}[Tuition Subsidies]
\label{exam:rct_fails}
A policy maker is interested in estimating the effect $\theta(p)$ of offering all high school graduates a
fixed subsidy of \$$p$ to attend college. To do so, they might consider running a small randomized controlled
trial: Given a small set of study participants, randomly assign half of them to receive a subsidy $p$ and
half of them not to, and then compare college enrollment rates among those two groups. As argued in
\citet{heckman1998general}, however, such an approach may badly over-estimate the effect of the subsidy
on enrollments because it fails to consider overall equilibrium effects on the college wage premium.

More formally, let $V(a, \, c)$ denote the average net value of enrolling in college,
where $a$ denotes the wage premium resulting from a college degree and $c$ the cost of attendance.
In general, we should expect $V$ to be monotonically increasing in $a$ and decreasing in $c$. The subsidy reduces
costs by $p$, and thus at first glance makes college more attractive. Where one needs to be careful, however, is in
recognizing that the college wage premium $a$ is not set in stone; rather, it is determined by labor market conditions.
If more people enroll in college, one may expect the labor market bargaining power of college graduates to
diminish, and for $a$ to decrease in response. Thus, if we believe the subsidy $p$ increases enrollments, we might
expect for $a(p)$ to be a (decreasing) function of $p$ due to equilibrium effects.

\begin{figure}
\begin{center}	
\includegraphics[width=0.7\textwidth]{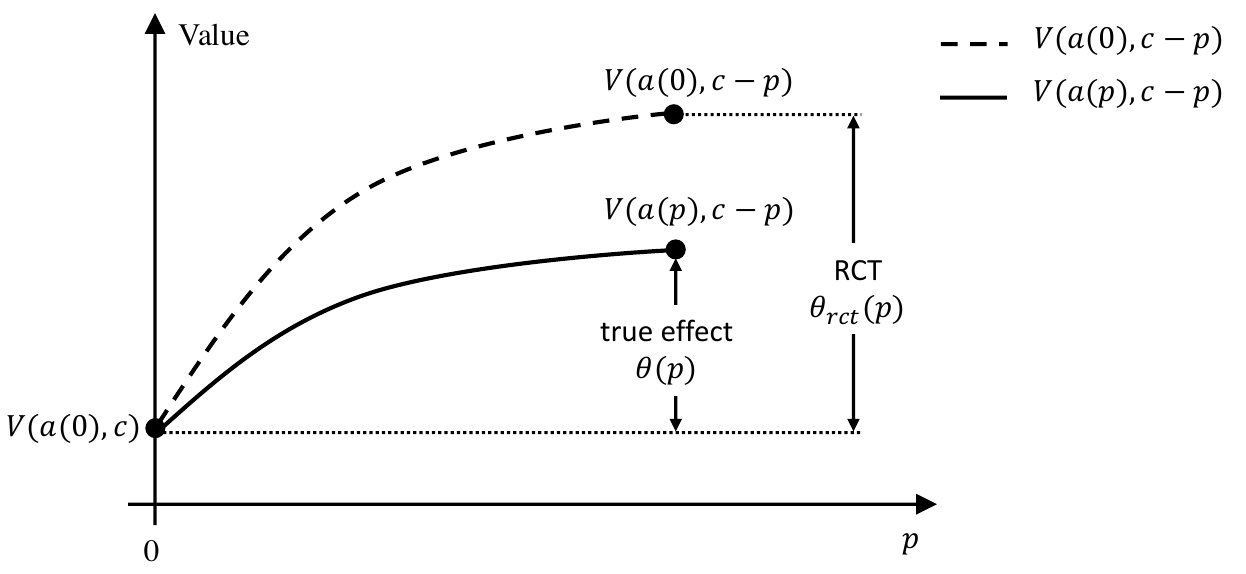}
\caption{Illustration of the failures of a randomized controlled trial under the presence of cross-unit interference in Example \ref{exam:rct_fails}.}
\label{fig:rct_fails}
\end{center}
\end{figure}

We are now ready to illustrate why a simple randomized trial falls short here (Figure \ref{fig:rct_fails}). The randomized trial only affects a
small number of study participants and does not capture changes in $a$; specifically, it measures
$\theta_{rct}(p) = V(a(0), \, c - p) - V(a(0), c)$. In contrast, the true effect of the subsidy should also
reflect its impact on the equilibrium college wage premium, i.e., $\theta(p) = V(a(p), \, c - p) - V(a(0),\,  c)$.
In general, for any subsidy $p > 0$, we should expect $a(p) < a(0)$ and so $V(a(0), \, c - p) > V(a(p), \, c - p)$, meaning
that the randomized trial will over-estimate the effect of the subsidy. On the quantitative front, \citet{heckman1998general}
discuss a setup where ignoring equilibrium effects would lead to estimates that are off by an order of magnitude.
\end{exam}

\subsection{Interference and Clustered Inference}

The question of how to run experiments in the presence of cross-unit interference has received considerable attention in the literature.
The simplest approach to dealing with interference is to assume that we can divide our experimental samples into disjoint clusters
that do not interfere with each other, and then to consider inference at the level of these clusters \citep{baird2018optimal,hudgens2008toward}.

One such example involves experimentation in internet ad auctions, where each auction consists of a keyword along with a set of advertisers who submit competing bids in order for their ads to be displayed when the keyword is queries by a user. There is cross-unit interference because the same advertiser or keyword may appear in multiple auctions. \citet{basse2016randomization} and \citet{ostrovsky2011reserve} make the observation that the auction type used for one keyword does not meaningfully affect how advertisers bid for other keywords. They then consider experiments that group auctions into clusters by their keywords and randomize auction formats across these keyword clusters, rather than across advertisers, as a means
to avoid problems with interference.
More broadly, in the context of tuition subsidies, this idea of cluster-level randomization could correspond to identifying communities that are relatively isolated
from each other and randomizing the interventions across communities rather than across individuals; or,
in the case of social networking, it could involve deploying different versions of a feature in different countries
and hope that the number of cross-border links is small enough to induce only negligible interference.

The limitation of such cluster-based approaches, however, is that the power of any experiment is limited by the number of non-interfering clusters available: For example, if a platform has 200 million customers in 100 countries, but chooses to randomize by country, then the largest effective sample size they can use for any experiment is 100, and not 200 million.
Recently, several authors have sought to improve on the power of such cluster-based approaches by considering methods
that allow interference to be captured by a generic graph, where two units are connected by an edge if the treatment assigned to one unit
may affect the other's outcome \citep{aronow2017estimating,athey2018exact,basse2019randomization,eckles2017design,leung2020treatment}.
Even in this general case, however, we typically need to assume that the interference graph is sparse, i.e.,
that most units do not interfere with each other. For example, \citet{leung2020treatment} assumes that the average degree of the
interference graph remains bounded.

\subsection{Accounting for Interference via Equilibrium Modeling}

In this paper, we propose an alternative approach to experimentation in stochastic systems, where a large number of, if not all, units interfere with one another.  For concreteness, we focus on the problem of setting supply side payments in a centralized marketplace, where available demand is randomly allocated to a set of available suppliers. In these systems, different suppliers interact via their effects on the overall supply-demand equilibrium: The more suppliers choose to participate in the marketplace, the less demand on average an individual supplier would be able to serve in equilibrium. The objective of the system designer is to identify the optimal payment that maximizes the platform's utility.  Note that conventional randomized experimentation schemes that assume no interference fail in this system: For example, if we double the per-transaction payments made to a random half of suppliers, these suppliers will be more inclined to participate and reduce the amount of demand available to the remaining suppliers, and thus reduce their incentives to participate. 

We consider a simple model of such a centralized marketplace, and
design a class of ``local'' experimentation schemes that---by carefully leveraging
the structure of the marketplace---enable us to optimize payments without
disturbing the overall market equilibrium. To do so, we perturb the
per-transaction payment $p_i$ available to the $i$-th supplier by a small
mean-zero shock, i.e., $p_i = p + \zeta \varepsilon_i$ where
$0 < \zeta \ll 1$ and $\varepsilon_i = \pm 1$ independently and uniformly at random.
A reduced form linear regression, one that estimates how the individual random shock $ \zeta \varepsilon_i$ affects supplier-$i$'s behavior, recovers
a certain marginal response function, which captures the supplier's sensitivity to payment changes against a fixed ambient market equilibrium. This marginal response, unfortunately, is not directly relevant for policy design, as it does not take into account the shift in market equilibrium should all suppliers receive the same payment change. However, in the limit where the number of suppliers is large, we show that a mean-field
model can be used to translate the output of this reduced form regression into an estimate of the gradient of
the platform's utility with respect to $p$. We can then use these gradient estimates to optimize $p$ via any
stochastic first-order optimization method, such as stochastic gradient descent and its extensions.

The driving insight behind our result is that, although there is dependence across the
behavior of a large number of units in the system, any such interference can only
be channeled through a small number of key statistics: In our example, this corresponds to the total
supply made available by all suppliers. Then, if we can intervene on individual units without
meaningfully affecting the key statistics, we can obtain useful information about the
system---at a cost that scales sub-linearly in the number of units.
The type of interference that we consider,  where the units experience global interference channeled through
a small number of key statistics, can manifest in a range of applications. We discuss some examples below. 

\begin{exam}[Ride Sharing]
\label{exam:2}
Ride sharing platforms match customers who request a ride with nearby freelance drivers who are both active
and not currently servicing another request. It is in the interest of the platform to have a reasonable amount
of capacity available at all times to ensure a reliable customer experience. To this end, the platform may seek
to increase capacity by increasing the rates paid to drivers for completing rides. And, when running experiments
on the rates needed to achieve specific capacity levels, the platform needs to account for interference.
If the platform in fact succeeds in increasing capacity by increasing rates---yet demand remains fixed---the expected utilization
of each driver will go down and so the drivers' expected revenue, i.e., the product of the rate and the expected utilization, will
not increase linearly in the rate. Thus, if drivers respond to expected revenue when choosing whether to work for a platform, as
empirical evidence suggests that they do \citep{hall2019pricing}, a platform that ignores interference effects will overestimate
the power of rate hikes to increase capacity.
However, as shown in our paper, we can accurately account for these interference effects via mean-field
modeling because they are all channeled through a simple statistic, in this case total capacity. 
\end{exam}

\begin{exam}[Congestion Pricing]
\label{exam:3}
A policy maker may want to identify the optimal toll for congestion pricing \citep[e.g.,][]{goh2002congestion}.
We assume that drivers get positive utility from completing a trip, but get negative utility both from congestion delays
and from paying tolls. Then, in studying the effect of a toll on congestion, the policy maker needs to address the
fact that drivers interfere with one another through the overall state of congestion on the road: If we raise the tolls
on a small subset of the drivers and hence discourage them from going on the road, those whose tolls remain
unchanged may experience  less congestion and hence be inclined to drive more. Therefore a policy maker that experiments
with a small sub-population, without taking into account interference effects, may obtain an overly
optimistic estimate of the true effect of a toll change when applied to all drivers. Again, however, all interference
is channeled through a single statistic---congestion---and so mean-field modeling can capture its effect. 
\end{exam}

\begin{exam}[Renewable Energy Subsidies]
\label{exam:4}
In an electricity whole sale market, energy producers (e.g., generators) and consumers
(e.g., utilities) make bids and offers in the day-ahead market, which is then
cleared in a manner that balances the aggregate regional supply and demand. The
operator of these markets, such as CAISO or ERCOT, may choose to provide subsidies or scheduling priorities
to encourage renewable generation \citep[see][]{isocal2009}. Suppose that the market
operator would like to know the effect of increasing subsidies on energy generation.
We expect that increased subsidies would increase both total and renewable energy
production; the question is by how much, and what the effect of interference will be.
It is plausible that the effect of subsidies on total supply will be mitigated by interference,
because increased production from one supplier will decrease demand available to others.
In contrast, interference may either mitigate or amplify the effect of subsidies on renewable
energy production: Amplification effects may occur if subsidies affect profitability in a way
that causes non-renewable producers to be replaced by new renewable entrants. In either
case, all interference effects are channeled through global capacity, and so can
be accounted for via mean-field modeling. 
\end{exam}

\subsection{Related Work}
\label{sec:rel_work}

The problem of experimental design under interference has received considerable attention in
the statistics literature. For example,
\citet{blake2014marketplace} document failures of the non-interference assumption due to an interaction
between treated and control customers in an experiment run by an online marketplace.
\citet{blundell2004evaluating} consider the effects of a job search program on employment outcomes,
and emphasize the importance of considering general equilibrium effects whereby
job offers given to program participants may substitute for job offers given to non-participants and
increased search activity from participants may lower equilibrium wages for less skilled individuals.
\citet{bottou2013counterfactual} describe difficulties in using randomized experiments to study internet
ad auctions: Advertisers participate in an auction to determine ad placements, and any intervention on one advertiser may change their behavior on the auction and thus affect the opportunities available to other advertisers.
In all these cases, simple randomized controlled trials would paint a misleading picture about the
effect of an overall policy change.

The dominant paradigm for working under interference has focused on robustness to
potential interference effects, and on defining estimands in settings where
some units may be exposed to spillovers from treating other units
\citep{aronow2017estimating,athey2018exact,baird2018optimal,basse2019randomization,eckles2017design,
hudgens2008toward,leung2020treatment,manski2013identification,sobel2006randomized,tchetgen2012causal}.
Depending on applications, the exposure patterns may be simple (e.g., the units are clustered
such that exposure effects are contained within clusters)
or more complicated (e.g., the units are connected in a network, and two units far from each other
in graph distance are not exposed to each others' treatments). Unlike this line
of work that seeks robustness to interference driven by potentially complex and unknown mechanisms,
the local randomization scheme proposed here crucially relies on having a stochastic
model that lets us explain interference. Then, because all inference
acts via a simple statistic, we can move beyond simply seeking
robustness to interference and can in fact accurately predict
interference effect using information gathered in equilibrium.

Another plausible approach would be to use structural estimation 
methods and directly estimate the whole underlying system, and subsequently use stochastic optimization to obtain the optimal decision. However, a full-blown structural estimation approach would be infeasible in our problem because it involves a large number of interacting units each with unknown features. In particular, as will be clear in Section \ref{sec:model}, we consider the interaction among a large number of units, and each unit's behavior depends on a random choice function drawn from a potentially large set of options. The set of problem parameters thus involves the shapes of every possible choice function, as well as sampling distribution with respect to which the function is drawn for each unit. Directly estimating these parameters can be very difficult, and as we show, is not needed if the final goal is to identify the optimal action. Instead, our approach will focus on estimating a small number of key statistics which turn out to be sufficient for performing optimization. Doing so allows us to side-step the scalability problem of the structural estimation approach and arrive at the optimal action in an efficient manner.

The idea that one can distill insights of a structural model down to the relationship
between a small number of observable statistics has a long tradition in
economics \citep[e.g.,][]{chetty2009sufficient,harberger1964measurement}.
This approach can often be used for practical counterfactual
analysis without needing to fit complicated structural models.
We are inspired by this approach, and here we use such an argument
for experimental design rather than to guide methods for observational study analysis.
At a high level, our paper also has a connection to results on learning in a setting where agents exhibit strategic behavior, 
including \citet{feng2018learning}, \citet{iyer2014mean}, and \citet{kanoria2017dynamic}, and in crowd-sourcing systems, including \citet{johari2017matching}, \cite{khetan2016achieving} and \citet{massoulie2018capacity}. 

Our approach to optimizing $p$ using gradients obtained from local experimentation intersects with the literature on
continuous-arm bandits (or noisy zeroth-order optimization), which aims to optimize
a function $f(x)$ by sequentially evaluating $f$ at points $x_1, x_2, \ldots$, and obtaining in return noisy versions
of the function values $f(x_1), f(x_2), \ldots$ \citep{bubeck2017kernel,spall2005introduction}. 
A number of bandit methods first generate noisy gradient estimates of the function by comparing adjacent function values, and subsequently use these estimates in a first-order optimization method \citep{flaxman2005online, ghadimi2013stochastic, jamieson2012query, kleinberg2005nearly, nesterov2017random}.  In our model, this approach would amount to estimating 
utility gradients via what we call global experimentation, i.e., by comparing the  empirical utilities observed at two different payment levels.
Compared to this literature, our paper exploits a cross-sectional structure not present in most existing zeroth-order models: We show that  our local experimentation approach, which offers slightly different payments across a large number of units, is far more efficient at estimating the gradient than global experimentation, which offers all units the same payment on a given day.
Such cross-sectional signals  would be lost if we abstracted away the multiplicity of units, and only treated the average payment as a decision variable to be optimized.
In Section \ref{sec:bandit}, we provide a formal comparison for the regret of a platform deploying our approach versus
a bandit-based algorithm, and establish sharp separation in terms of rates of convergence.

The limiting regime that we use, one in which the system size tends to infinity, is often known as the mean-field limit. It has a long history in the study of large-scale stochastic systems, such as the many-server regime in queueing networks \citep{bramson2012asymptotic, halfin1981heavy, stolyar2015pull, tsitsiklis2012power, vvedenskaya1996queueing} and interacting particle systems \citep{graham1994chaos, mezard1987spin, sznitman1991topics}.
A key property of this mean-field limit is that, while changes to the behavior of a single unit may have significant impact on other units in a finite system, such interference diminishes as the system size grows and, in the limit, the behaviors among any finite set of units become asymptotically independent from one another, a phenomenon known as the propagation of chaos \citep{bramson2012asymptotic, graham1994chaos, sznitman1991topics}. This asymptotic independence property underpins the effectiveness of our local experimentation scheme, and ensures that small, symmetric payment perturbations do not drastically alter the equilibrium demand-supply dynamics.

Mean-field-inspired approaches have also been used in game theory to analyze equilibria in the presence of a large number of players by assuming that the agents respond to a certain average behavior of the system \citep{adlakha2015equilibria,hopenhayn1992entry,jovanovic1988anonymous, weintraub2008markov}; the equilibrium notion we use also falls under this category. In contrast to the existing literature, the main focus of our work lies in using mean-field limits to drive learning and experimentation.

\section{Designing Experiments under Equilibrium Effects}
\label{sec:motivation}

For concreteness, we focus our discussion on a simple setting inspired by a centralized marketplace
for freelance labor that operates over a
number of periods. In each period, the high-level objective of the decision maker (i.e., operator of
the platform) is to match demand with a pool of potential suppliers in such a manner
that maximizes the platform's expected utility. To do so, the decision maker
offers payments to each potential supplier individually, who in turn decides whether
to become active/available based  upon their belief of future revenue. Our main
question is how the decision maker can use  experimentation to efficiently
discover their revenue-maximizing payment, despite not knowing the detailed
parameterization of the model, and the presence of substantial stochastic uncertainty.

We formally describe a flexible stochastic model in Section \ref{sec:model}; here, we briefly
outline a simple variant of our model that lets us highlight some key properties of our approach.
Each day $t = 1, \, ..., \, T$ there are $i = 1, \, ..., \, n$
potential suppliers, and demand for $D_t$ identical tasks to be accomplished. A central platform chooses
a distribution $\pi_t$, and then offers each supplier random payments \smash{$P_{it} \simiid \pi_t$} they
commit to pay for each unit of demand served. The suppliers
observe both $\pi_t$ a state variable $A_t$ that can be used to accurately anticipate demand $D_t$
(e.g., $A_t$ could capture local weather or events);
however, the platform does not have access to $A_t$. Given their knowledge of $P_{it}$ and $A_t$,
each supplier independently chooses to become ``active''; we write $Z_{it} = 1$ for active suppliers
and $Z_{it} = 0$ else. Then, demand $D_t$ is randomly allocated to active suppliers.

Our key assumption is that each supplier chooses to become active based on
their expected revenue conditionally on being active, and furthermore that they
do so via stationary reasoning \citep{hopenhayn1992entry}. Each supplier first computes $q_{A_t}(\pi_t)$, their expected
allocation rate (rate at which they will be matched with demand) conditionally on being active 
and given $A_t$ and $\pi_t$. They then decide whether to become active by comparing the
expected revenue $P_{it} q_{A_t}(\pi_t)$ with a random outside option. We refer to this as a
stationary model of supplier choice as we implicitly assume that suppliers don't take into account the
effect of their own decision to become active on their expected allocation rate.
This is often taken to be a reasonable assumption in large stochastic systems
\citep{adlakha2015equilibria,chetty2009sufficient,weintraub2008markov}.

\begin{figure}
\begin{center}
\includegraphics[width=0.48\textwidth]{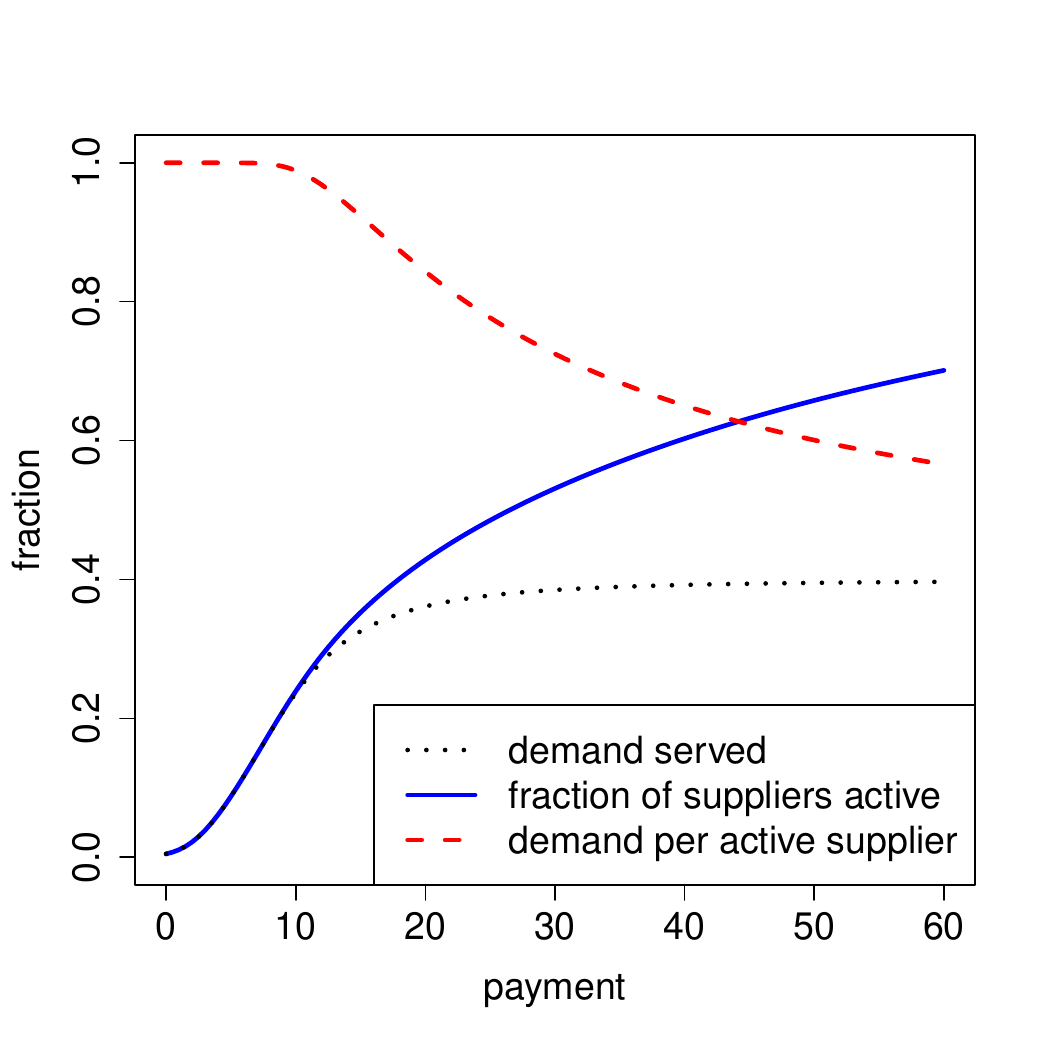}
\caption{Example of large-sample behavior of market, conditionally on a realization of $A$.
We show $\mu_A(p)$, the fraction of suppliers that choose to become active,
 $q_A(p)$ the expected amount of demand served per active supplier and $\mu_A(p) q_A(p)$ the expected
amount of demand served (expressed as a multiple of the maximum capacity that would be available
if all suppliers were active). The example is simulated in the mean-field limit, i.e., with number of potential
suppliers $n$ growing to infinity, such that $\EE{D/n \cond A} = 0.4$. Individual supplier preferences
are logistic \eqref{eq:logistic} with $\alpha = 1$ with outside option $\log(B_i/20) \sim \nn\p{0, \, 1}$.
Supply-demand matching is characterized via the allocation function \eqref{eq:queue} with $L = 8$, visualized in Figure \ref{fig:raf}.}
\label{fig:example}
\end{center}
\end{figure}

The form of $q(\cdot)$ depends on both the amount of available supply and demand, and the
efficiency with which supply can be matched with demand; see Section \ref{sec:model} for
an example based on a queuing network. Finally, the platform's utility $U_t$ is given by the revenue
from the demand served minus payments made to suppliers. 
Figure \ref{fig:example} shows a simple example of an equilibrium resulting from this
model in the limit as $n$ gets large in a setting where all suppliers are offered the same payment $p$, 
for a specific realization of demand $D$. We see that, as $p$ gets larger, the active supply gets larger
than demand and the utilization of active suppliers goes down.

Conversely, our assumption that the platform cannot observe the daily state variable $A_t$
is made to ensure that our learning problem is robust and performs well even if the state variables are unavailable or difficult to estimate accurately.
In practice, of course, it is plausible that a platform may have access to partial---but not full---information
about $A_t$. Here, we focus on the statistically most difficult setting where the
platform is oblivious to $A_t$ and thus can only learn how to set $p$ via experimentation,
as this setting enables us to establish a crisp separation between different approaches to learning
and to highlight our core methodological contributions.
However, all methods considered here can be adapted to leverage partial information about $A_t$,
and further work that investigates how best to leverage such information would be of considerable interest.

Before presenting our proposed approach to learning $p$ below, we first briefly review why standard
approaches fall short. The core difficulty in our model comes from the interplay between network
effects and market-wide demand fluctuations induced by the $A_t$.

The network effects break what one might call classical A/B-experimentation. Suppose that,
on each day $t = 1, \, ..., \, T$, the platform chooses a small random fraction of suppliers
and offers them an experimental payment $p_{\text{exp}}$, while everyone else gets offered the status
quo payment $p_{\text{default}}$. We could then try to use the behavior of suppliers offered $p_{exp}$ to estimate
expected profit at $p_{\text{exp}}$, and then update the default payment. This approach allows for cheap
experimentation because most of the suppliers get offered $p_{\text{default}}$. However, it will not
consistently recover the optimal payment because it ignores feedback effects: When we raise payments, more
suppliers opt to join the market and so the rate at which any given supplier is matched with demand goes
down---and this attenuates the payment-sensitivity of supply relative to what is predicted by A/B testing.

Conversely, the market-wide demand fluctuations due to $A_t$ degrade global optimization schemes
that use payment variation across days for learning; such algorithms are equivalent to continuous-armed
bandit algorithms considered in the optimization literature \citep{spall2005introduction}.
Suppose that, on each day $t = 1, \, ..., \, T$, we randomly chose a payment $p_t$ and made it
available to all suppliers, and then observed realized profits $U_t$. We could then try
estimate profit gradients by comparing $U_{t}$ to $U_{t-1}$. The problem is that, due to variation
in daily context, the variation in per-supplier profit $U_t/n$ given the chosen payment $p_t$ is always
of constant order, even in very large markets (i.e., in the limit $n \rightarrow \infty$); for example,
in a ride-sharing setting, if day $t-1$ is rainy and day $t$ is sunny, then the effect of this weather
change on profit may overwhelm the effect of any payment change deployed by the
platform.\footnote{Of course, the platform may try to correct for contexts, e.g., by matching days
with similar values of $A_t$ with each other. One currently popular way of doing so in the technology
industry is using synthetic controls \citep{abadie2010synthetic}. In practice, however, this approach
may be difficult to implement, and will remain intractably noisy unless the platform can observe
the full context $A_t$ and use it to essentially perfectly predict demand. As discussed above, our goal
in this paper is to develop methods for learning that are driven purely by experimentation, and that do not rely
on the platform being able to accurately observe $A_t$.}
The upshot is that the platform cannot learn anything via global experimentation unless it considers
large changes to the payments $p_t$ that it offers to everyone. And such wide-spread payment changes
are impractical for several reasons: They are expensive, and difficult to deploy.

\subsection{Local Experimentation}

Our goal is to use high-level information about the stochastic system described above to design
a new experimental framework that lets us avoid the problems of both approaches described above:
We want our experimental scheme to be consistent for the optimal payment (like global experimentation), but also
to be cost-effective (like classical A/B testing) in that it only requires small perturbations to the status quo.

The driving insight behind our approach is that it is possible to learn about the relationship between profit
and payment via unobstrusive randomization by randomly
perturbing the payments $P_{it}$ offered to supplier $i$ in time period $t$. We propose
setting
\begin{equation}
\label{eq:local_randomization}
P_{it} = p_{t} + \zeta \varepsilon_{it}, \ \ \ \ \varepsilon_{it} \simiid \cb{\pm 1}
\end{equation}
uniformly at random, where $\zeta > 0$ is a (small) constant that governs the magnitude of the perturbations,
and regressing market participation $Z_{it}$ on the payment perturbations $\varepsilon_{it}$.
This regression lets us recover the \emph{marginal response function}, i.e.,
the average payment sensitivity of a supplier in a situation where only they get different payments but
others do not; see Section \ref{sec:MRF} for a formal definition.

This marginal response function is not directly of interest for optimizing $p$, as it ignores feedback
effects. However, we find that---in our setting---this quantity
captures relevant information for optimizing payments. More specifically we show in Section \ref{sec:MRF}
that, provided we have good enough understanding of system dynamics to be able to anticipate
match rates given the amount of supply and demand present in the market, in the mean-field limit where the market size grows, we can use
consistent estimates of the marginal response function to derive consistent estimates of
the actual payment-sensitivity of supply that accounts for network effects.
Furthermore, we show in Section \ref{sec:learning} that this approach enables us to optimize
payments using vanishingly small-scale experimentation as the market gets large (i.e., we can
take $\zeta$ in \eqref{eq:local_randomization} to be very small when $n$ is large).

\begin{figure}
\begin{center}
\begin{tabular}{cc}
\includegraphics[width=0.48\textwidth]{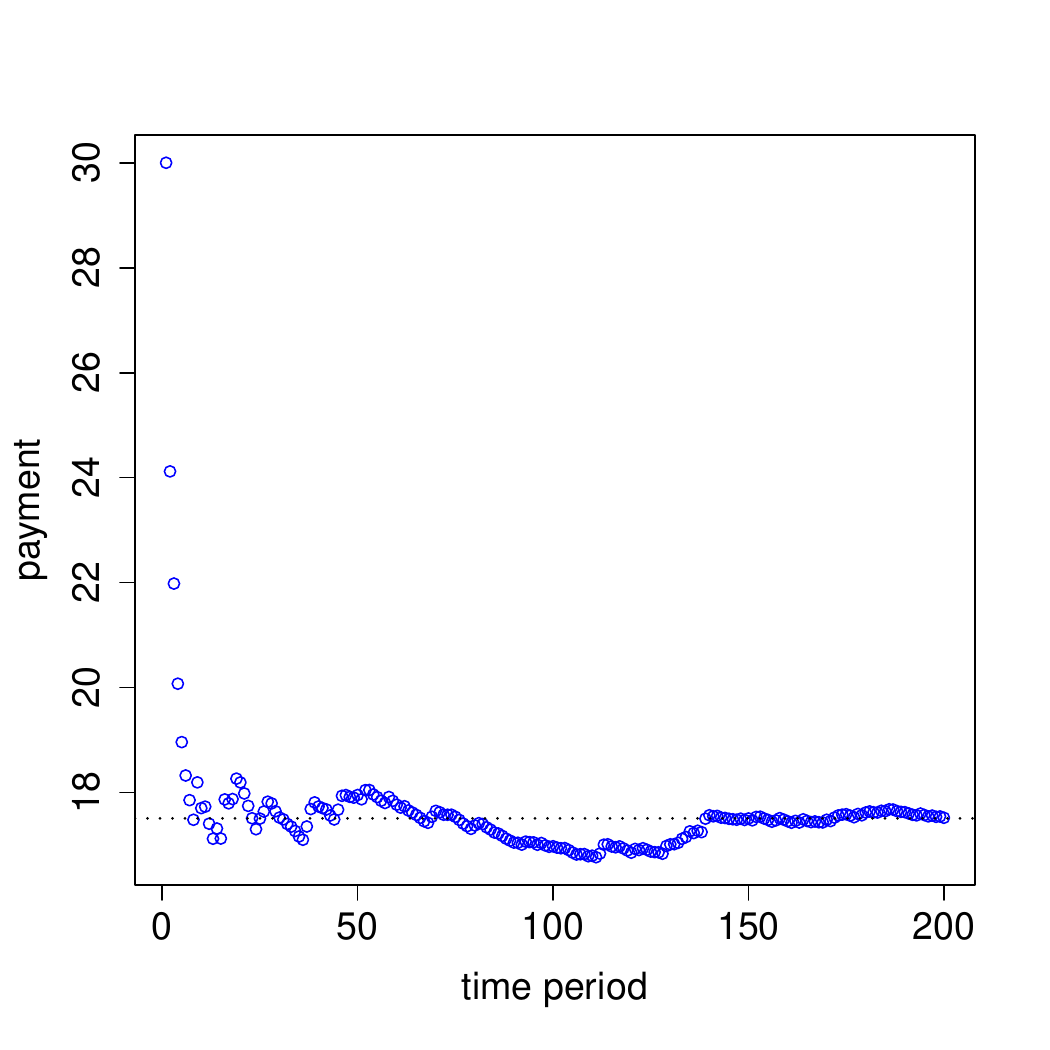} &
\includegraphics[width=0.48\textwidth]{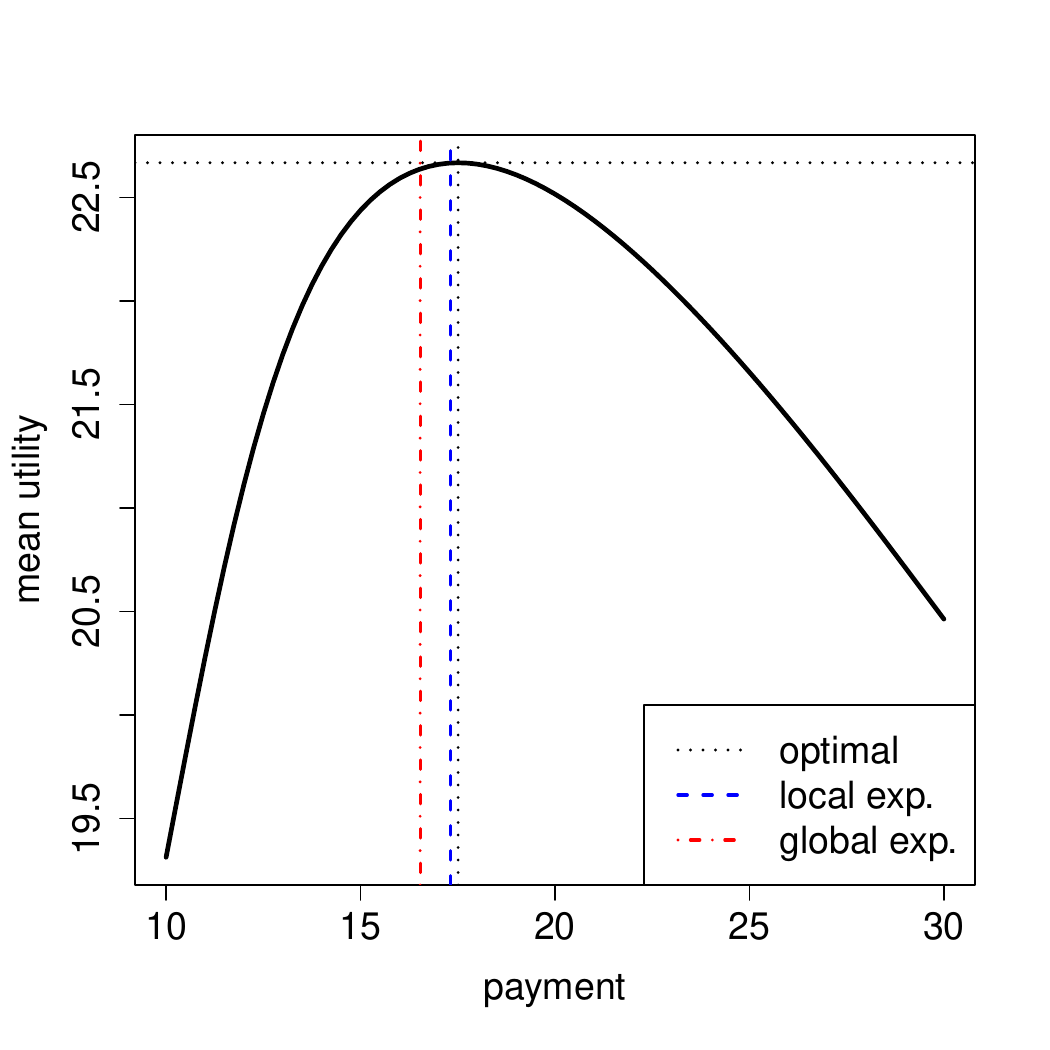}
\end{tabular}
\end{center}
\caption{Results from learning $p$ via local experimentation. The worker preference functions are as
in Figure \ref{fig:example}; the daily contexts are such that $\EE{D/n \cond A} \sim \text{beta}(15, 35)$.
The platform utility function is linear as in Lemma \ref{lemm:concave_u}, with $\gamma = 100$.
We learned gradients based on local randomization \eqref{eq:local_randomization} with $\zeta = 0.5$, and
then optimized payments via gradient descent as in \eqref{eq:OMD} with a step size $\eta = 20$ and $I = (-\infty, \, \infty)$.
The left panel shows the convergence of the $p_t$ to the value $p^*$ that optimizes mean utility.
The right panel compares the average value of $p_t$ over the last 100 steps of our algorithm
to both a payment $\hat{p}$ learned via global experimentation and the optimal payment $p^*$.}
\label{fig:sgd}
\end{figure}

Figure \ref{fig:sgd} shows results from our local experimentation approach
on a simple simulation experiment in the setting of Figure \ref{fig:example}, where
the scaled demand $\EE{D/n\cond A}$ follows a $\text{beta}(15, \, 35)$ distribution. We initialize the
system at $p_1 = 30$, and then each day run payment perturbations as in \eqref{eq:local_randomization}
to guide a payment update using an update rule described in Section \ref{sec:algo}.
We see that the system quickly converges to a near-optimal payment of around 17.

\begin{figure}
\begin{center}
\begin{tabular}{cc}
\includegraphics[width=0.48\textwidth]{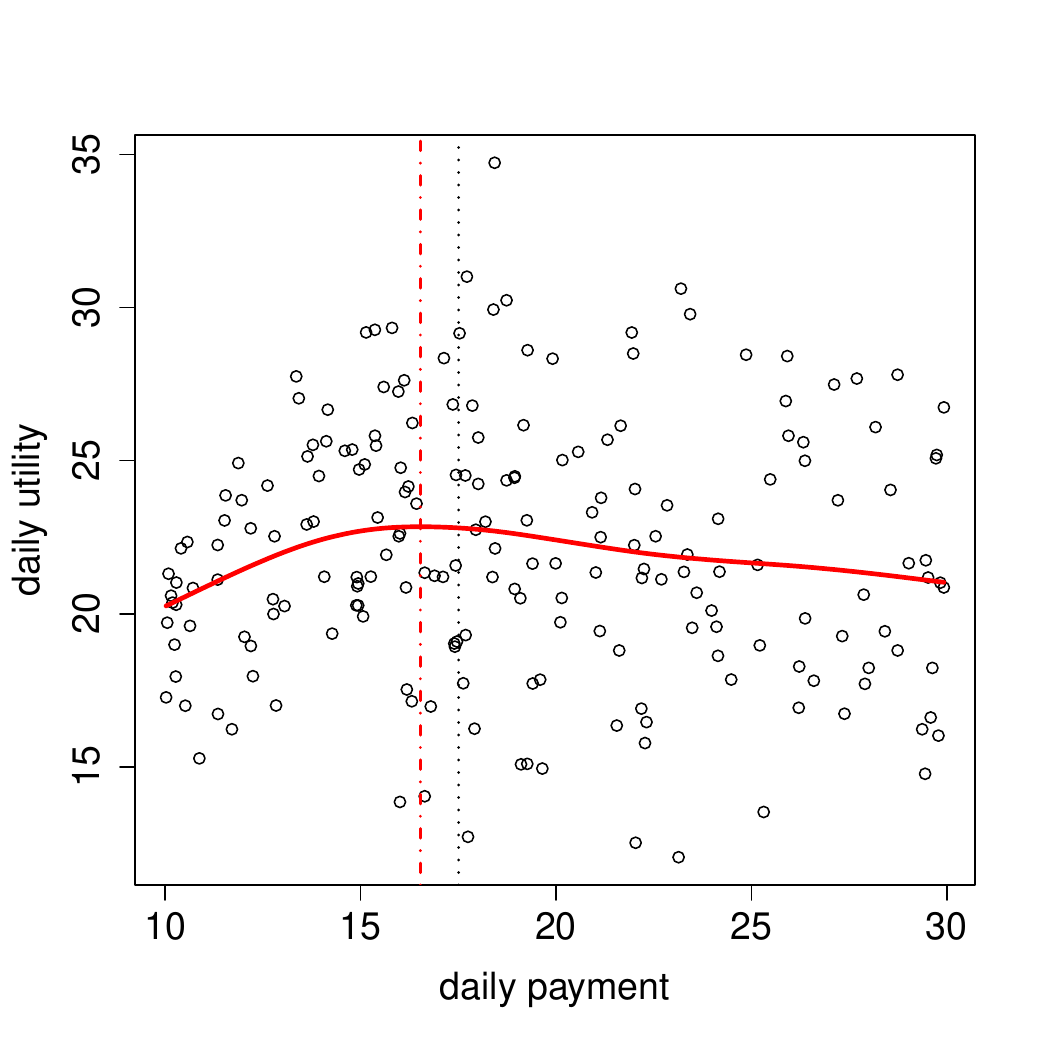} &
\includegraphics[width=0.48\textwidth]{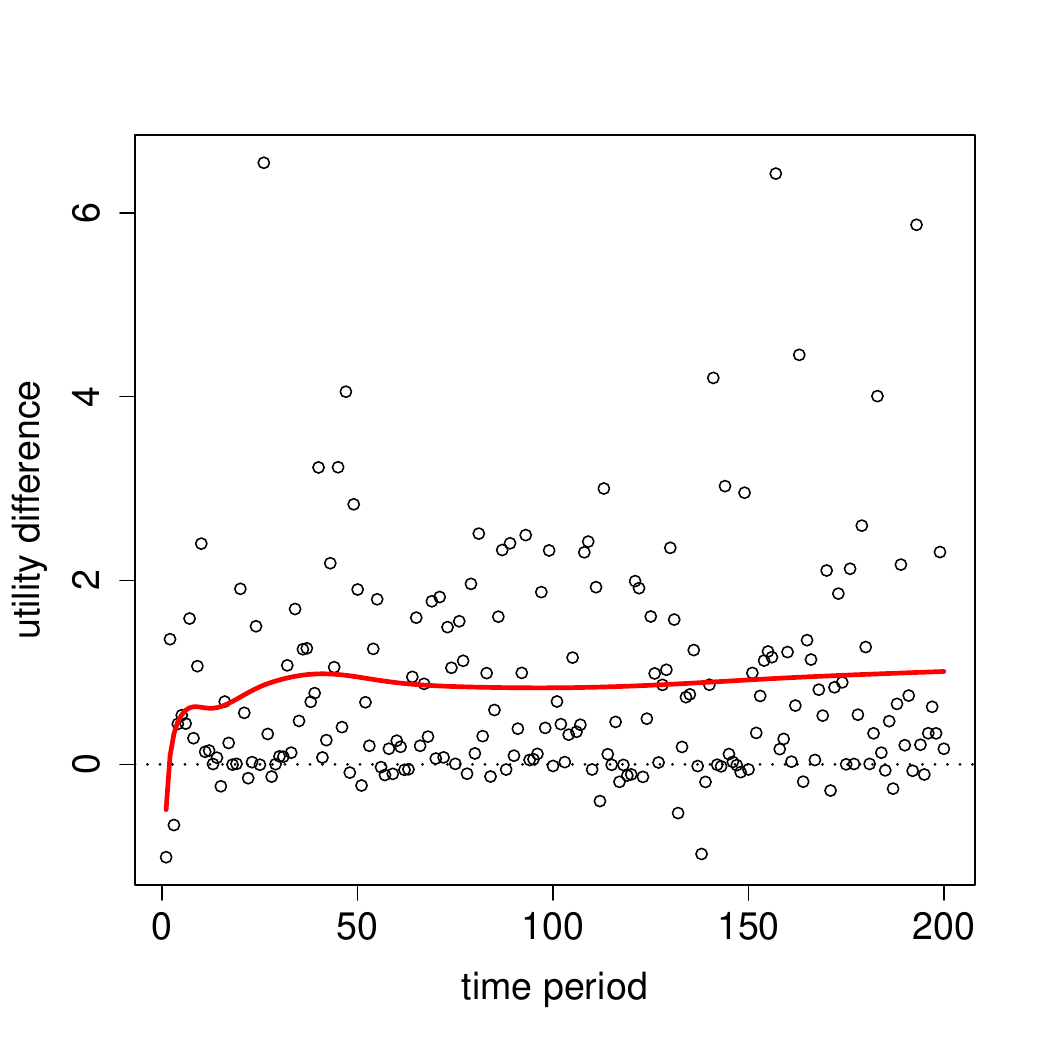}
\end{tabular}
\end{center}
\caption{Results from learning $p$ via global experimentation. The left panel shows pairs $(p_t, \, U_t)$
resulting from daily experiments, along with both the resulting $\hat{p}$ (dash-dotted line) and the optimal
$p^*$ (dotted line). The right panel shows the
(scaled) difference in daily utility between our local experimentation approach and the global experimentation
baseline (both approaches worked using the same demand sequence $D_t$).} 
\label{fig:baseline}
\end{figure}

We also compare our results to what one could obtain using the baseline of global experimentation,
where we randomize the payment $p_t \sim \text{Uniform}(10, \, 30)$ in each time period and measure resulting platform utility $U_t$, and then choose the final
payment $\hat{p}$ by maximizing a smooth estimate of the expectation of $U_t$ given $p_t$. The left panel of
Figure \ref{fig:baseline} shows the resulting $(p_t, \, U_t)$ pairs, as well as the resulting $\hat{p}$.
As seen in the right panel of Figure \ref{fig:sgd}, the final $\hat{p}$ obtained via this method is a reasonable estimate
of the optimal $p$.

The major difference between the local and global randomization schemes is in the resulting cost of experimentation.
In Section \ref{sec:cost} we show that our local experimentation scheme pays a vanishing cost for randomization; the
only regret relative to deploying the optimal $p$ from the start is due to the rate of convergence of gradient descent.
In contrast,the cost of experimentation incurred for finding $\hat{p}$ via global experimentation is huge, because it
needs to sometimes deploy very poor choices of $p_t$ in order to learn anything. And, as shown in the right panel of
Figure \ref{fig:baseline}, after the first few days, the global experimentation approach in fact
systematically achieves lower daily utilities $U_t$ than local experimentation.
In Section \ref{sec:simu} we consider further numerical comparisons of local and global experimentation, as well as
variants of global exploration that balance exploration and exploitation to improve in-sample regret.

\begin{rema}[Relationship to Batched Bandits]
\label{rema:batch}
Our model bears resemblance to batched multi-arm bandits \citep{esfandiari2019batched,gao2019batched,perchet2016batched} and batched online optimization \citep{bubeck2019complexity,duchi2018minimax}, where an analyst sequentially picks multiple arms to pull for one batch at a time. In particular, administering an intervention to a unit in our model could be seen as analogous to pulling one arm in batched bandits, or sampling an unknown function at a particular point in batched online optimization. There is, however, a fundamental distinction between our model and the predominant model for batched bandits. Existing work on batched bandits does not allow for interference within batches: The action assigned to one unit in a batch does not directly affect the outcome observed for another unit in the batch. In contrast, the presence of cross-unit interference within batches (or, for us, within days) is at the heart of our model: The outcome of a unit not only depends on their own intervention, but also on the interventions experienced by other units on the same day. Thus, existing results on batched bandits and online optimization cannot be used to reason about how best to deploy heterogeneous incentives to different suppliers in order to converge to a good choice of $p$ in our setting. 
\end{rema}

\section{Model: Stochastic Market with Centralized Pricing}
\label{sec:model}

We now present the general stochastic model we use to motivate our approach.
All random variables are assumed to be independent across the periods and, within each period, are independent from one another unless otherwise stated. We will consider a sequence of systems, indexed by $n\in \N$, where in the $n$-th system there are $n$ potential suppliers. We will refer to $n$ as the market size.
All variables in our model are thus implicitly dependent on the index, $n$, which we
denote using the superscript $(n)$, e.g., $q\nth$.
We sometimes suppress this notation when the context is clear.  In the rest of the section, we will focus on describing the model in a single time period. 

\paragraph{Demand}
To reflect the reality that demand fluctuations may not concentrate with $n$, we allow for
a random stochastic global state $A$ drawn from a finite  set $\calA$. The global state affects demand,  and is known to market participants (suppliers), but not to the platform (or the platform cannot react to it). For example, in a ride sharing example, $A$ could capture the effect of weather (rain / shine) or  major events (conference, sports game, etc.).
Conditionally on the global state $A = a$, we assume that  demand, $D$, is drawn from distribution $D \sim F_a$. We further assume that the demand scales proportionally with respect to the market size $n$, and that it concentrates after re-scaling by $1/n$.
In particular, we assume that there exists $\{d_a\}_{a\in \calA} \subset \rp$, such that for all $a\in \calA$, $\EE{D/n | A=a} = d_a$ for all $n\in \N$, 
\begin{equation}
\lim_{n \to \infty}\EE{\p{D/n - d_a}^2 \bbar A=a} = 0, 
\label{eq:Dconv0}
\end{equation}
and
\begin{equation}
\pb\p{ D/n \notin  [d_a/2, 2d_a] \bbar A=a} = o(1/n), 
\label{eq:Dconv}
\end{equation}
and as $n\to \infty$. In general, we  will use the sub-script $a$ to denote the conditioning that the global state $A=a$.

\paragraph{Matching Demand with Suppliers}
Depending on the realization of demand, all or a subset of the suppliers will be selected to serve the demand. In particular, the matching between the  potential suppliers and demand occurs in three rounds:

\begin{enumerate}[]
\item  {\it Round 1:} The platform chooses a {payment distribution}, $\pi$, and draws payments
$P_i \simiid \pi$  for $i = 1, 2, \ldots, n$. Then, for each supplier $i$, the platform announces both
the payment $P_i$ and the underlying distribution $\pi$, with the understanding that the supplier will be
compensated with $P_i$ for every unit of demand that they will be matched with eventually.

\item {\it Round 2:}  Suppliers choose whether to they want to be {active}.   A supplier will not be matched with any demand if they choose to be inactive.  We write $Z_i \in \cb{0, \,1}$ to denote whether the $i$-th participant chooses to participate in the marketplace, and write $T = \sum_{i = 1}^n Z_i$ as the total number of active suppliers.  The mechanism through which a supplier determines whether or not to become active will be described shortly. 
 
\item
 {\it Round 3:}  The platform employs some mechanism that randomly matches demand with active suppliers.
\end{enumerate}
Denote by $S_i$ the amount of demand that an active supplier $i$ will be able to serve, and define 
\begin{equation}
 \Omega(d, t) \bydef \E [ S_i  \bbar D=d, T=t], 
 \label{eq:alloc_prelim}
 \end{equation} 
as the expected demand allocation to an active supplier under the payment distribution $\pi$, conditional on the total demand being $d$ and total active suppliers being $t$. 
We allow for a range of possible matching mechanisms, but assume that in the limiting regime where $t$ and $d$ are large, $\Omega(d,t)$
converges to a ``regular allocation function'' that only depends on the ratio between the demand and active suppliers, $d/t$.

 \begin{defi}[Regular Allocation Function] A function $\omega: \rp \to [0,1]$ is a regular allocation function if it satisfies the following: 
 \label{def:RAF}
 \begin{enumerate}
\item $\omega(\cdot)$ is smooth, concave and non-decreasing.
 \item $\lim_{x\to 0} \omega(x) = 0$ and $\lim_{x\to \infty}\omega(x) \leq 1 $.
 \item $\lim_{x\to 0}\omega'(x) \leq 1$.
 \end{enumerate}
 \end{defi}

The condition of $\omega$ being concave corresponds to the assumption that the marginal difficulty with which additional demand can be matched does not decrease as demand increases. The condition that $\lim_{x\to \infty}\omega(x) \leq 1$ asserts that the maximum capacity of all active suppliers be bounded after normalization.

\begin{assu} 
\label{ass:Omega_tou}
The function  $\Omega: \rp^2  \to \rp$ satisfies the following: 
\begin{enumerate}
\item $\Omega(d,t)$ is non-decreasing in $d$, and non-increasing in $t$. 
\item There  exists a bounded error function $l: \rp^2 \to \rp$ with 
\begin{equation}
\label{eq:Omegaconvg}
\abs{l(d,t)} = o\p{1/\sqrt{t} + 1/\sqrt{d}}, 
\end{equation}
such that $\Omega(d,t)= \omega(d/t) + l(d,t)$ for all $t,d \in \rp$, where $\omega(\cdot)$ is a regular allocation function. 
\end{enumerate}
\end{assu}

We provide below an example system in which the allocation rates are given by a regular allocation function (Definition \ref{def:RAF}). 

\begin{figure}
\centering
\includegraphics[width=0.44\textwidth]{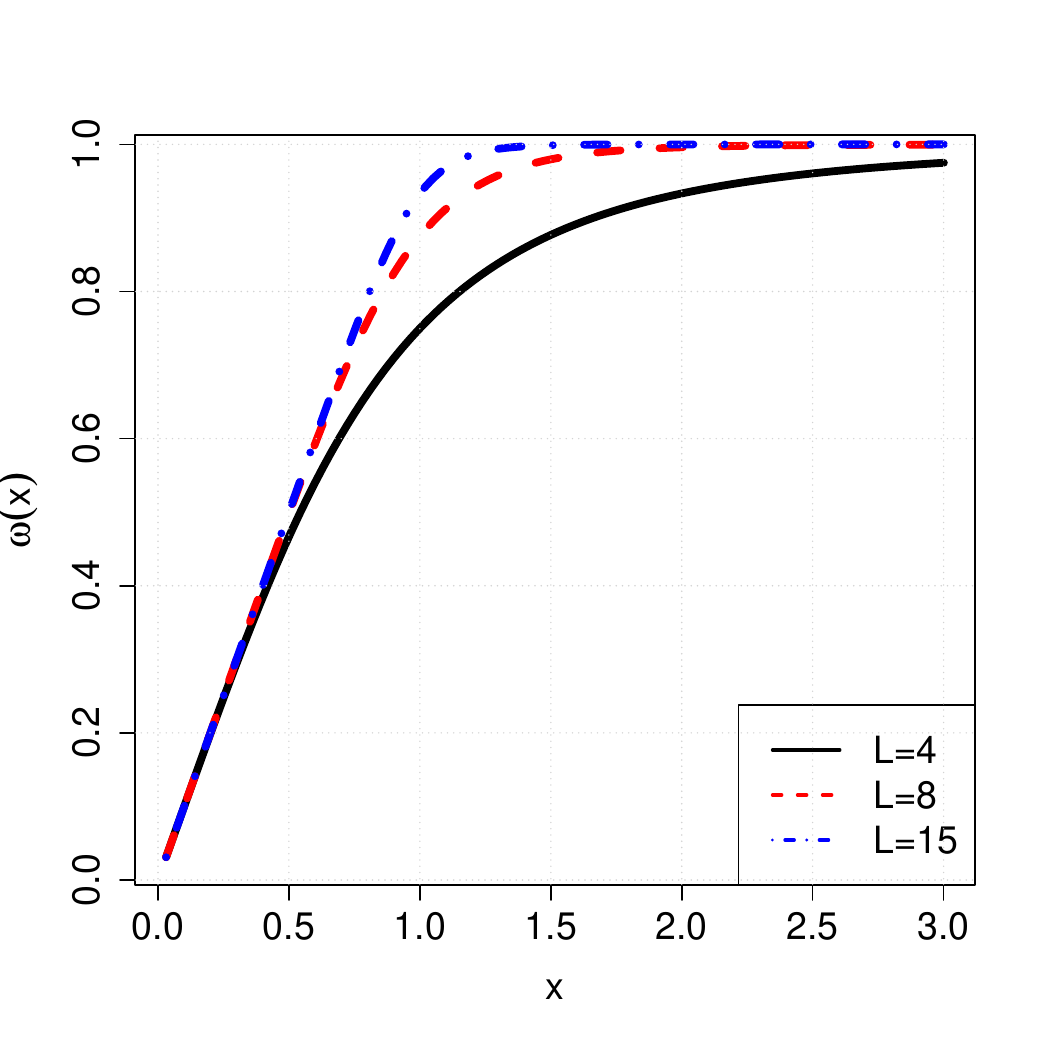}
\caption{Examples of the regular allocation function $\omega(\cdot)$ in Example \ref{exp:queues} under different values of capacity $L$. }
\label{fig:raf}
\end{figure}

\begin{exam}[Regular Allocation Function Example: Parallel Finite-Capacity Queues] 
\label{exp:queues}
Consider a service system where each active supplier operates as a single-server $M/M/1$ queue with a finite capacity, $L \in \N$, $L\geq 2$. A request that arrives at a queue is accepted if and only if the queue length is less than or equal to $L$, and is otherwise dropped. We assume that all servers operate at unit-rate, so that a request's service time is an independent exponential random variable with mean 1. Each unit demand generates an independent stream of requests which is modeled by a unit-rate Poisson process, so that the aggregate arrival process of requests  is Poisson  with rate $D$  (by the merging property of independent Poisson processes). When a new request is generated within the system, the platform routes it to one of the $T$ queues selected uniformly at random. The random routing corresponds, for instance, to a scenario where both the incoming requests and active suppliers are scattered across a geographical area, and as such, requests are assigned to the nearest server. 

Within this model, each active supplier effectively functions an $M/M/1$ queue with service rate $1$ and arrival rate $D/T$. Due to the capacity limit at $L$, some requests may be dropped if they are assigned to a queue currently at capacity. Using the theory of $M/M/1$ queues, it is not difficult to show that  \citep[e.g.,][eq.~5.6]{spencer2014queuing} if we denote $D/T$ by $x$, then the rate at which requests are processed by a server, corresponding to the allocation rate,  is given by
\begin{equation}
\label{eq:queue}
\begin{split}
\omega(x) = \left\{
                \begin{array}{ll}
                  \frac{x-x^L}{1-x^L},&  \quad x \neq 1, \\
                  1-\frac{1}{L}, &\quad  x = 1.
                \end{array}
              \right.
\end{split}
\end{equation}
Numerical examples of $\omega(\cdot)$ are given in Figure \ref{fig:raf}.  Note that $\omega(\cdot)$ satisfies all conditions in Definition \ref{def:RAF} and is hence a regular allocation function. Finally, we may generalize the model to where the suppliers are partitioned into $k$ equal-sized groups, so that each sever operates at speed $Tkm$. The corresponding allocation function would have the same qualitative behavior. 
\end{exam}

\paragraph{Supplier Choice Behavior}
We assume that each supplier takes into account their expected revenue in equilibrium when making the decision of whether or not to become active. In particular, the  of supplier $i$ becoming active is given as follows, where $T $ is the equilibrium number of active suppliers:
\begin{equation}
\mu\nth_a(\pi) \bydef \PP[\pi]{Z_i=1 \bbar A=a} = \EE[\pi]{f_{B_i}( P_i \, \EE[\pi]{\Omega(D, T) \bbar A=a}) \cond A=a}.
\label{eq:choice}
\end{equation} 
Here, $\EE[\pi]{\Omega(D, T)\, |\, A=a}$ is the expected amount of demand served by each supplier given the platform's choice of $\pi$, 
and thus $P_i \,  \EE[\pi]{\Omega(D, T) \, |\, A=a}$ is the expected revenue of the $i$-th supplier in equilibrium.\footnote{For now, assume that such equilibrium distribution is well defined, and we will justify its meaning rigorously in a moment.}
Note that the choice model \eqref{eq:choice} is stationary in that each supplier only considers the average behavior of
other marketplace participants when choosing whether or not to enter. In particular, suppliers do not consider the effect of their own
entry decision on the system, or combinatorial interactions between other marketplace participants. Similar types of stationary assumptions, also known as mean-field or oblivious equilibrium, are common in game theoretic models involving a large number of players where each player's influence on the overall system dynamic is vanishingly small \citep{hopenhayn1992entry, weintraub2008markov}, and can be formally justified by showing how stationary equilibrium converges to the true Nash equilibrium in the limit as the system size tends to infinity \citep{adlakha2015equilibria}. 

Here, $B_i$ is a private feature that captures the heterogeneity across potential suppliers, such as a supplier's cost, or noise in their estimate of the expected revenue. We assume that the $B_i$'s  are drawn i.i.d.~from a set $\calB$ whose distribution may depend on $A$. The choice function  $f_b(x)$ represents the  of the supplier becoming active, when their private feature is $b$ and expected equilibrium revenue is $x$. We assume the family of choice functions $\{f_b(\cdot)\}_{b\in \calB}$ satisfy certain regularity properties detailed below.

\begin{assu}
\label{ass:choicef}
We assume that supplier choices are determined by the stationary choice model \eqref{eq:choice}. Furthermore,
for all $b\in \calB$, we assume that the choice function $f_b(\cdot)$ takes values in $[0,1]$, is monotonically
non-decreasing, and twice differentiable with a uniformly bounded second derivative. 
\end{assu}

\noindent
Below is one example of a family of choice functions that satisfies Assumption \ref{ass:choicef}: 

\begin{exam}[Logistic Choice Function]
A popular model in choice theory is the logit model (cf.~Chapter 3 of \cite{train2009discrete}), which, in our context, corresponds to the choice function being the logistic function: 
\begin{equation}
\label{eq:logistic}
\PP{Z_i=1 \cond P_i, \, \pi, \, A} = \frac{1}{1 + e^{-\alpha\p{P_i \EE[\pi]{\Omega(D, T) \bbar A} - B_i}}},
\end{equation}
where $\alpha>0$ is a parameter and the private feature $B_i$ takes values in $\rp$ and represent the break-even cost threshold of supplier $i$. In this example, the supplier's decision on whether to activate will depend on whether their expected revenue exceeds their break-even cost. The sensitivity of such dependence is modeled by the parameter $\alpha$. Note that in the limit as $\alpha \to \infty$, the probability of the event $Z_i=1$ conditionally on $P_i$, $\pi$ and $A$ is either 0 or 1. That is, a supplier will choose to be active if and only if they believe their expected revenue from Round 2 will exceed the break-even threshold $B_i$. 
\end{exam}

\paragraph{Platform Utility and Objective}  The platform's utility is defined to be the difference between revenue and total payment: 
\begin{equation}
U = R(D, \, T) - \sum_{i = 1}^n P_i Z_i S_i, 
\label{eq:reward}
\end{equation}
where $S_i$ is the amount of demand that a supplier would serve if they become active, and
$R(D, T)$ is the platform's expected  revenue, with equilibrium active supply size $T$ and total demand $D$.  Analogously to the case of $\Omega(D,T)$, we will assume that the revenue function $R$  is approximately linear in the sense that,  for some function $r$,  $R( D, \, T) \approx r(D/ T)T$ when $T$ and $D$ are large. More precisely, assume the following:

\begin{assu} 
\label{ass:Rtou}
There  exists a bounded error function $l: \rp^2 \to \rp$ with $\abs{l(d,t)} =  \smash{o(1/\sqrt{t} + 1/\sqrt{d})}$ such that 
\begin{equation}
R(d,t)  = \p{r(d/t) - l(d,t)}t, \quad \mbox{for all $t,d \in \rp$}, 
\label{eq:Rconvg}
\end{equation}
where $r: \rp \to \rp$ is a smooth function with bounded derivatives.
\end{assu}
As an example, the platform could receive a fixed amount $\gamma$ from each unit of demand served, in which case we have $R(D,T) = \gamma (T \Omega(D,T))$. 
Given this notation, we write the platform's expected utility in the $n$-th system as 
\begin{equation}
\label{eq:chia}
\begin{split}
u_a\nth (\pi)&= \frac{1}{n} \EE[n]{U \cond A = a}, \quad \mbox{and} \quad u\nth(\pi) = \EE[n]{u_A\nth(\pi)}.
\end{split}
\end{equation}
Denote by $\delta_x$ the Dirac measure with unit mass on $x$. We consider
two different objectives for the decision maker (i.e., platform operator). First, they may want to control
\emph{regret}, and deploy a sequence of payment distributions $\pi$ whose utility nearly matches that of
the optimal \emph{fixed payment}, $p^*$. Second, they may want to estimate $p^*$.
In Section \ref{sec:learning}, we provide results with guarantees along both objectives.

\paragraph{Symmetric Payment Perturbation}
An important family of payment distributions that will be used repeatedly throughout the paper is that of symmetric payment perturbation. Let $\{\varepsilon_i\}_{i \in \N}$ be a sequence of i.i.d.~Bernoulli random variables with  $\pb(\varepsilon_i = -1)= \pb(\varepsilon_i = +1)= \frac{1}{2}. $ Fix $p> \zeta>0$. We say the payments  are $\zeta$-perturbed from $p$, if
\begin{equation}
P_i = p+ \zeta\varepsilon_i , \quad i \in \N. 
\label{eq:zEta}
\end{equation}
In what follows, we will use $\pi_{p,\zeta}$ to denote the payment distribution when payments are $\zeta$-perturbed from $p$,
\smash{$\mu\nth_a(p, \zeta)$} to denote \smash{$\mu\nth_a(\pi_{p,\zeta})$}. The meanings of \smash{$\mu\nth_a(p, \zeta)$},
\smash{$u\nth_a(p,\zeta)$}, etc., are to be understood analogously. When $\zeta=0$, we may omit the dependence on
$\zeta$ and write, for instance, \smash{$\mu\nth_a(p)$} in place of \smash{$\mu\nth_a(p,0)$} or \smash{$\mu\nth_a(\pi_{p,0})$}.

\begin{rema}[What does the platform know?] 
\label{rema:who_knows}
Our model assumes that the platform has detailed knowledge of the allocation
mechanics, but cannot anticipate the behaviors of market participants that drive of supply and demand.
More specifically, we assume that the platform knows
the regular allocation function $\omega$ (Definition \ref{def:RAF}), and its  pre-limit version, $\Omega$ \eqref{eq:alloc_prelim};
the limiting platform utility function $r$ (Assumption \ref{ass:Rtou}), and its pre-limt version, $R$ \eqref{eq:reward};
as well as the payment scheme it chooses to use, i.e., $p$, $\zeta$, and the realizations of the random  perturbations, $\{\varepsilon_i\}_{i=1, \ldots, n}$.
However, the platform cannot anticipate the global state $A$, the demand $D$, or the distribution of supplier
choice functions $f_{B_i}(\cdot)$; rather, all it can do is collect after-the-fact measurements of
$D$ and $\{Z_i\}_{i=1, \ldots, n}$, the set of active suppliers.
Finally, we implicitly assume that the the global state $A_t$ has no effect on the system beyond their assigned time period, and that the platform knows this fact.

This modeling choice reflects an understanding that it is realistic for a platform to have a
good handle on the mechanics of the marketplace it controls, but it is implausible for it
to have an in-depth understanding of the beliefs and preferences of all marketplace participants. For
example, in the case of ride sharing, it is plausible that a platform could get good at modeling congestion,
but less plausible that the platform could fully understand and anticipate how all its drivers may respond
to various policy changes.

The fact that we take the platform to be completely oblivious to the global state $A_t$ puts us in
an extreme setting, where the platform's learning must be purely driven by randomization in $p$. We chose this extreme setting primarily for two reasons. First, it crystallizes the difficulty of the learning problem, and highlights the value of local experimentation relative
to global experimentation baselines. Second, a platform's knowledge of $A_t$, if any, is likely to be noisy and inaccurate, and
it is often difficult for a platform to learn efficiently in practice by matching historical data using noisy estimates of their corresponding contexts.
Therefore, it is of considerable practical value to devise a robust learning algorithm that works well without relying on the platform's ability to infer the global state $A_t$. 

In practice, of course, the platform may have some information about the global state $A$; for example,
we may assume that the platform observes a set of covariates $X$ that capture some aspects of $A$
(e.g., we could have $X = \Xi(A)$ for some lossy function $\Xi$). In such a setting, the information $X$
could be used for variance reduction and/or learning better policies that exploit heterogeneity explained by $X$.
It would be of considerable interest to study a covariate-enriched variant of our approach that
allows the platform to use such information to learn better policies; however,
we leave this line of investigation to follow-up work.
\end{rema}

\subsection{Mean-Field Asymptotics}
\label{sec:limits}

The stochastic model described above in general admits complex dynamics that are not amenable to exact analysis. Fortunately, we show in this sub-section that in the mean-field limit where the number of suppliers is large, various key equilibrium quantities converge to tractable objects described by a mean-field model. To start, we first provide a formal definition of the equilibrium active supply size, $T$, and verify existence and uniqueness.

\begin{defi}[Active Supply Size in Equilibrium]
We say that a random variable $T$ is an \emph{equilibrium supply size}, if, when all suppliers make activation choices according to \eqref{eq:choice}, the resulting distribution for the  number of active suppliers equals that of $T$. 
\end{defi}

\begin{lemm}
\label{lem:sa_unique} 
Suppose that the conditions in Assumptions \ref{ass:Omega_tou}, \ref{ass:choicef} and \ref{ass:Rtou} hold.  Fix $p>0$, $\zeta \in [0, p)$, and $a \in \calA$. Let the payment distribution $\pi$ be defined on $\rp$. Then, conditional on $A=a$, the equilibrium active supply size exists, is unique, and follows a Binomial distribution.
\end{lemm}

Next, we define some quantities that will play a key role in our analysis, and verify that they converge
to tractable mean-field limits. The first quantity we consider is the equilibrium number of active suppliers
$\mu_a^{(n)}(p)$, as defined in \eqref{eq:choice}. Second, we define the function $q(\cdot)$, which captures
the expected amount of demand matched to each supplier if the total number of suppliers were exogenously drawn
as a binomial $(n, \, \mu)$ random variable rather than determined by the equilibrium:
\begin{equation}
\label{eq:qdef}
q\nth_a(\mu)   =  \EE{ \Omega \p{D, X}|A=a}, \ \ \ \ X \sim \text{Binomial}(n, \, \mu).
\end{equation}

\begin{lemm}
\label{lem:conv_fluid}
Under the conditions of Lemma \ref{lem:sa_unique}, for all $a\in \calA$, and $p,\mu \in \rp$, the following hold:
\begin{align}
\label{eq:muconvg}
&\limn \mu_a \nth (p) =  \mu_a(p), \\
\label{eq:qconvg}
&\limn  q\nth_a(\mu) =  \omega(d_a/\mu), \\
\label{eq:uconvg}
&\limn u\nth_a(p) =  u_a(p) =\p{ r(d_a/\mu_a(p)) - p\omega(d_a/\mu_a(p))}\mu_a(p), \\
\label{eq:qprimeconvg}
&\limn (q\nth_a)' (\mu) =  - \omega'(d_a/\mu) \frac{d_a}{\mu^2}, 
\end{align}
where $\omega(\cdot)$ and $r(\cdot)$ are described in Definition \ref{def:RAF} and Assumption \ref{ass:Rtou}, respectively.  In \eqref{eq:muconvg}, the limit $\mu_a(p)$ is the
only solution to $\mu = \EE{f_{B_1}\p{p \omega(d_a/\mu)} \cond A=a}$.
\end{lemm}

Finally, the following result, proven in Appendix \ref{sec:concave_u}, establishes conditions under which the limiting utility functions $u_a(p)$
are concave, thus enabling us to globally optimize utility via first-order methods.

\begin{lemm}
\label{lemm:concave_u}
Let $f_a(\cdot)$ be the average choice function:  $f_a(x) = \EE{f_{B_1}(x) \cond A=a}. $  Fix $\gamma>0$, $c_0 \in (0, \gamma)$ and $a\in \calA$. Suppose the following holds: 
\begin{enumerate}
\item We have a linear revenue function, $ r(x) = \gamma \omega(x)$ . 
 \item Let $\underline{x} = \inf_{p \in (c_0, \gamma) } pq_a(\mu_a(p))$ and $\overline{x} = \sup_{p \in (c_0, \gamma) } pq_a(\mu_a(p))$.  The  average choice function  $f_a(\cdot)$ satisfies 
\begin{enumerate}
\item $f_a(\cdot)$ is strongly concave in the domain $(\underline{x}, \overline{x})$. 
\item $f_a(\underline{x}) - f'_a(\underline{x})\underline{x} \geq 0$, or, equivalently, that there exists a differentiable, non-negative concave function $\tilde{f}(\cdot)$, such that $\tilde{f}(\underline{x}) = f_a(\underline{x})$  and $\tilde{f}'(\underline{x})\leq f'_a(\underline{x})$. 
\end{enumerate} 
\item The allocation function $\omega(\cdot)$ is strongly concave in the domain $\p{d_a/\mu_a(c_0), d_a/\mu_a(\gamma)}$. 
\end{enumerate}
Then, under the conditions of Lemma \ref{lem:sa_unique},  the limiting platform utility $u_a(\cdot)$ is strongly concave in the domain $ (c_0,\gamma)$. 
\end{lemm}

\subsection{The Marginal Response Function}
\label{sec:MRF}

Finally, as discussed in Section \ref{sec:motivation}, a key quantity that motivates our approach
to experimentation is the marginal response function, $\Delta(p)$, which captures the average payment sensitivity
of a supplier in a situation where only they get different payments but others do not
(meaning that there are no network effects).

\begin{defi}[Marginal Response Function] Fix $n\in \N$, $a\in \calA$ and $p >0$. The marginal response function is defined by
\begin{equation}
\label{eq:Delta}
\Delta\nth_a(p)  =   \qatp \EE {f'_{B_1}\p{ p q\nth_a\p{\mu\nth_a(p)) }}  \cond A=a}.
\end{equation}
\end{defi}

This marginal response function $\Delta$ plays a key role in our analysis for the following reasons.
First, as shown in the following section,  in the mean-field limit as $n\to \infty$, $\Delta$ is easy to estimate using small random payment
perturbations that do not meaningfully affect the overall equilibrium.
Second, provided we have a good enough understanding of the underlying system dynamics
to know the appropriate allocation function $\omega(\cdot)$, we can use consistent estimates
of $\Delta$ to estimate the true payment sensitivity of supply that accounts for feedback effects,
$d\mu(p)/dp$. 
This fact is formalized in the following result. We note that, other than $\Delta$, all terms on the right-hand side of \eqref{eq:mup2}
are readily estimated from observed data by taking averages. 

\begin{lemm} 
\label{lem:limit_mup}
Under the conditions of Lemma \ref{lem:sa_unique}, for any $a\in \calA$ and $p \in \rp$, we have that
\begin{align}
\label{eq:mup1}
\frac{d}{dp}\mu_a\nth(p)  =  \frac{\Delta\nth_a(p)}{1-p \Delta\nth_a(p)  \, q^{(n)'}_a\p{\mu_a\nth(p)}/\qatp  } \text{ for any $n \geq 1$}.
\end{align}
Furthermore, this relationship carries through in the mean-field limit, 
\begin{align}
\label{eq:deltaconvg}
&\limn \Delta\nth_a (p) =  \Delta_a(p) \bydef  \omega\p{d_a/\mu_a(p)} \EE{f'_{B_1}\p{p\omega\p{d_a/\mu_a(p)}}  \cond A=a}, \\
\label{eq:mup2}
& \lim_{n\to \infty} \frac{d}{dp}\mu\nth_a(p) =  \mu_a'(p) = 
 \Delta_a(p) \,\Big/\, \p{1+\frac{pd_a \Delta_a(p)\omega'\p{d_a/\mu_a(p)}}{\mu_a(p)^2 \omega(d_a/\mu_a(p))}}.
\end{align}
\end{lemm}

In addition to powering our approach to experimentation, the result of Lemma \ref{lem:limit_mup} also provides
qualitative insights about the drivers of interference in our model.
If there were no interference among the suppliers, then the gradient \smash{$(d/dp)\mu_a(p) $}
would have coincided with the marginal response \smash{$\Delta_a(p)$}; but due to interference, the
gradient is attenuated by an interference factor $1 + R_a(p)$, where
\begin{equation}
R_a(p) = \Sigma_a^\Delta(p)\Sigma_a^\Omega(p), \ \  \underbrace{\Sigma_a^\Delta(p) = \frac{p \Delta_a(p)}{\mu_a(p)},}_{\text{scaled marginal sensitivity}} \ \   \underbrace{\Sigma_a^\Omega(p) = \frac{d_a}{\mu_a(p)} \frac{\omega'(d_a/\mu_a(p))}{\omega(d_a/\mu_a(p))}.}_{\text{scaled matching elasticity}}
\end{equation}
We thus observe the following:
\begin{itemize}
\item The interference factor is negligible when the ``scaled marginal sensitivity'' \smash{$\Sigma_a^\Delta(p)$} is small, i.e.,
the marginal response function is small relative to the current supply $\mu_a(p)$. Note that $p \Delta_a(p)$ is a scale-free version
our marginal response function that is invariant to rescaling $p$.
\item The interference factor is negligible when the ``scaled matching elasticity'' \smash{$\Sigma_a^\Omega(p)$} is small, i.e.,
the elasticity of the matching function $\omega(\cdot)$ is small relative to the current ratio of supply to demand $\mu_a(p) / d_a$.
In particular, because $\omega(\cdot)$ is concave and bounded by assumption,
we can verify that  \smash{$\Sigma_a^\Omega(p)$} is small whenever demand far exceeds supply, i.e. $d_a / \mu_a(p) \gg 1$;
see Proposition \ref{prop:weakconcave} stated below and proven in Appendix \ref{app:prop:weakconcave}.
\item The interference factor is non-negligible when neither of the above conditions hold.
\end{itemize}
These observations are aligned with what one might have anticipated based on qualitative arguments.
For example, interference effects clearly cannot matter
if marketplace participants are overall unresponsive to changes in $p$, and this is exactly what we found in the first bullet point.
Meanwhile, one might have expected for the effect of interference to be more pronounced when there is more intense
competition among the suppliers than when there is enough demand to keep all suppliers busy, and this conjecture is
well in line with our finding in the second bullet point.

\begin{prop}
\label{prop:weakconcave}
Let $g: \rp \to \rp$ be concave with piece-wise continuous derivative $g'$.
Then $xg'(x)  \leq g(x)$ for all $x > 0$. If moreover $0 < \lim_{x \rightarrow \infty} g(x) < \infty$,
then $\lim_{x \rightarrow \infty} x g'(x) / g(x) = 0$.
\end{prop}

\section{Learning via Local Experimentation}
\label{sec:learning}

We present our main results in this section. The main framework we adopt for learning payments is based on first-order optimization. 
First, we show in Section \ref{sec:gradient} that our local experimentation approach enables us to construct an asymptotically accurate estimate of the utility gradient at a given payment $p$, in the mean-field   limit as $n \rightarrow \infty$. Then, we use these gradient estimates to update the payment using a form of gradient ascent, and show that their performance is superior to what can be achieved via classical continuous-armed bandit and zeroth-order optimization algorithms. Specifically, we establish in Section \ref{sec:algo} an $O(1/T)$ upper bound for the rate of convergence to the optimal platform utility under our algorithm. In Section \ref{sec:cost}, we study the cost of the local experimentation needed to estimate utility gradients, and verify that it scales sub-linearly in $n$. Finally in Section \ref{sec:bandit}, 
we compare our results to those available to classical continuous-armed bandits, and show that it is not possible to achieve the $O(1/T)$ convergence rate within the classical bandit framework. 
Throughout this section, we focus on optimizing utility in the mean-field limit, while verifying that finite-$n$ errors have an asymptotically vanishing effect on learning.

\subsection{Estimating Utility Gradients}
\label{sec:gradient}

Recall that, in our model, there are two sources of randomness. First, there is the stochastic global context
$A \in \acal$, which affects overall demand. In the context of ride-sharing, $A$ could capture multiplicative
demand fluctuations due to weather or holidays. Second, there is randomness due to decisions
of individual market participants. This second source of error decays with market size size $n$.
Our goal here is to verify that local experimentation allows us to eliminate  errors of the second
type via concentration as the market size $n$ gets large. Conversely, because the context $A$ affects everyone in the
same way, there is no way to average out the effect of $A$ without collecting data across many days. 

Define  $\bZ = \frac{1}{n}\sum_{i=1}^n Z_i$ and $\bD = D/n$.  As discussed in Section \ref{sec:motivation} our proposal starts for perturbing individual
payments as in \eqref{eq:zEta}, and then estimating the regression coefficient \smash{$\hDelta$} of
market participation $Z_i$ on the perturbation $\zeta_n \varepsilon_i$, i.e.,
\begin{equation}
\hDelta = \zeta_n^{-1} \sum_{i = 1}^n (Z_i - \bZ)(\varepsilon_i - \bar{\varepsilon}) \,\big/\, \sum_{i = 1}^n(\varepsilon_i - \bar{\varepsilon})^2.
\label{eq:margResEst}
\end{equation}
Our first result below relates this quantity $\hDelta$ we can estimate via local randomization to
a quantity that is more directly relevant to estimating payments, namely the payments derivative of $u$
conditionally on the global state $A$.

\begin{theo}
\label{thm:main}
Suppose the conditions of Lemma \ref{lem:sa_unique} hold.   Let 
\begin{equation}
 \hUps =  \hDelta/\lt(1+ \frac{p \bD \hDelta \omega'\p{\bD/ \bZ}}{\bZ^2 \omega\p{\bD / \bZ}}\rt), 
\end{equation}
and
\begin{equation}
\hGamma = \hUps \left[ r\p{\bD/\bZ} - p \omega\p{\bD/\bZ} - \p{ r'\p{\bD/\bZ} - p \omega'\p{\bD/\bZ}} \bD/\bZ\right]  - \omega\p{\bD/\bZ}\bZ.
\label{eq:Gamma_def}
\end{equation}
Then, assuming that the perturbations scale as $\zeta_n = \zeta n^{-\alpha}$ for some $0 < \alpha <  0.5$,
\begin{equation}
\label{eq:Gamma}
\limn \PP{\abs{ \hGamma - \frac{d}{dp}u_A\p{p}} > \varepsilon} = 0,
\end{equation}
for any $\varepsilon > 0$.
\end{theo}

\begin{rema}[Population-wide Experimentation \& Symmetric Perturbation]
\label{rema:symPertrub}

It is instructive to note that our experimentation scheme in \eqref{eq:zEta} has two distinguishing features that depart from a conventional approach to A/B testing that would gave a small subset of suppliers an $\varepsilon$ increase in payment while keeping payments in the rest of the population unchanged. First, our perturbation is symmetric across the units (zero-mean perturbation), whereas in classical settings those in a treatment group may receive asymmetric, and possibly identical, treatments.
Second, we experiment across the entire population as opposed to a small sub-population.

These features are in fact deliberate and interdependent, and the rationales are as follows.
The perturbations being symmetric ensures that our experimentation scheme does not meaningfully shift the overall supply-demand equilibrium ($\mu\nth_a(p) $), which in turn allows us to circumvent the impact of cross-unit interference. 
Moreover, as we show in Section \ref{sec:cost}, the symmetric perturbations  lead to a small cost of experimentation: Roughly speaking, the effect of paying half of the population $\varepsilon$ more is roughly neutralized by simultaneously paying the other half $\varepsilon$ less.
Meanwhile, the fact that we experiment on the whole population enables us to attain reasonable power using small enough
perturbations $\varepsilon$ such as not to be biased by the curvature of the supplier-specific choice functions $f_b(\cdot)$.

If perturbing the price carries a large fixed cost and the analyst wishes to apply local experimentation only to a small set of the population, one could also consider letting the random shocks $\epsilon_i$ be equal to $0$ with non-zero probability, and restrict the rest of the analysis on the sub-population who has received a non-zero shock.
This would also amount to a valid local experimentation scheme. However, we note that the power of our approach to estimate the
marginal response function depends on $\Var{\epsilon_i}$ and so, in order the maintain a given level of power (i.e., to ensure that
$\Var{\epsilon_i}=1$), the platform would need to use larger shocks $\varepsilon_{i}$ for those suppliers that receive non-zero shocks.
This in turn may expose the analyst to larger approximation errors from linearly approximating a curved function.
\end{rema}

\subsection{A First-Order Algorithm}
\label{sec:algo}

Our key use of Theorem \ref{thm:main} involves optimizing for a utility-maximizing $p$. At every time period $t$,
$\hGamma_t$ is a consistent estimate of the gradient of $u_{A_t}(\cdot)$ at $p_{t-1}$, and we can plug it into
any first-order optimization method that allows for noisy gradients. The proposal below is a variant
of mirror descent that allows us to constraint the $p_t$ to an interval $I$ \citep[e.g.,][]{beck2003mirror}.
We need to specify a step size $\eta$, an interval $I = [c_-, \, c_+]$, and an initial payment $p_1$.
Then, at time period $t = 1, \, 2, \, ...$,
we do the following:
\begin{enumerate}
\item Deploy randomized payment perturbations \eqref{eq:zEta} around $p_{t}$ to estimate $\hGamma_t$ as in \eqref{eq:Gamma_def},
\item Perform a gradient update\footnote{Note that, without the constraint to the interval $I$, this update is
equivalent to basic gradient descent with $p_{t+1} = p_t + 2\eta \hGamma_t / (t+1)$.}
\begin{equation}
\label{eq:OMD}
p_{t+1} = \argmin_p\cb{\frac{1}{2\eta} \sum_{s = 1}^t s (p - p_{s})^2 - \theta_t p : p \in I}, \ \ \ \ \theta_t = \sum_{s = 1}^t s \hGamma_s.
\end{equation}
\end{enumerate}
The following result shows that if we run our method
for $T$ time periods in a large
marketplace and the reward functions $u_{a}(\cdot)$ are strongly concave,
then the utility derived by our first-order optimization scheme is competitive with any
fixed payment level $p$, up to regret that decays as $1/t$.\footnote{In \eqref{eq:consistency},
we up-weight the regret terms $u_{A_t}(p) -  u_{A_t}(p_t)$ in later time periods to emphasize
their $1/t$ rate of decay. One could also use an analogous proof to verify that the unweighted
average regret is bounded on the order of
$T^{-1} \sum_{t = 1}^T \p{u_{A_t}(p) -  u_{A_t}(p_t)} = \oo_P(\log(T))$.}

\begin{theo}
\label{thm:learning}
Under the conditions of Theorem \ref{thm:main}, suppose we run the above learning algorithm
for $T$ time periods and that $u_{a}(\cdot)$ is $\sigma$-strongly concave over the interval $p \in I$ for all $a$.
Suppose, moreover, that we run \eqref{eq:OMD} with step size $\eta > \sigma^{-1}$ and that the gradients
of $u$ are bounded, i.e., $\abs{u'_{a}(p)} < M$ for all $p \in I$ and $a \in \acal$. Then
\begin{equation}
\label{eq:consistency}
\limn \PP{\frac{1}{T} \sum_{t = 1}^T t \p{u_{A_t}(p) -  u_{A_t}(p_t)}\leq \frac{\eta M^2}{2}} = 1,
\end{equation}
for any $p \in I$ and $T \geq 1$.
\end{theo}

The above result doesn't make any distributional assumptions on the contexts
$A_t$; rather, \eqref{eq:consistency} bounds the regret of our payment sequence $p_1, \, p_2, \, ...$
along the realized sample path of $A_t$ relative to any fixed oracle. We believe this aspect of our result
to be valuable in many situations: For example, if $A_t$ needs to capture weather phenomena that have a
big effect on demand, it is helpful not to need to model the distribution of $A_t$, as the weather may
have complex dependence in time as well as long-term patterns. However, if we are willing to
assume that the $A_t$ are independent and identically distributed, Theorem \ref{thm:learning}
also implies that an appropriate average of our learned payments is consistent for the optimal payment
via online-to-batch conversion \citep{cesa2004generalization}.

\begin{coro}
\label{coro:IID}
Under the conditions of Theorem \ref{thm:learning}, suppose moreover that the $A_t$
are independent and identically distributed and let $u(p) = \EE{u_{A_t}(p)}$.
Then, for any $\delta > 0$,
\begin{equation}
\label{eq:IID}
\limsup_{n \rightarrow \infty} \ \PP{ (p^* - \bar{p}_T)^2
\leq \frac{\eta M^2}{\sigma T}  \p{16 \log\p{\delta^{-1}} + 4}} \geq 1 - \delta,
\end{equation}
where $p^* = \argmax\cb{u(p) : p \in I}$ and $\bar{p}_T = \frac{2}{T(T+1)} \sum_{t = 1}^T t \, p_t$.
\end{coro}

\subsection{The Cost of Experimentation}
\label{sec:cost}

Our argument so far has proceeded in two parts. In Section \ref{sec:gradient} we showed we could consistently
use local experimentation to estimate gradients of the utility function $u_a(p)$. Then, in Section \ref{sec:algo}, we gave
bounds on the regret that updates payments $p_t$ via gradient descent---as though the platform could observe gradients
$u_a'(p)$ at no additional cost. Here, we complete the picture, and show that local experimentation in fact induces
negligible excess cost as we approach the mean-field limit. In general, a platform that randomizes payments around $p_t$
will make lower profits than one that just pays everyone $p_t$;\footnote{This is because randomization will not affect active
supply size to first order, but suppliers randomized to higher payments are more likely to be active. Randomization thus increases
the average per-unit payment the platform needs to give to suppliers without increasing the amount of demand the platform is able to serve.}
the result below, however, shows that this excess cost decays quadratically in the magnitude of payment perturbations $\zeta$.

\begin{theo}
\label{thm:cost}
Under the conditions of Theorem \ref{thm:main} there are constants $C, \, \alpha > 0$ such that
\begin{equation}
\label{eq:cost}
\frac{1}{T}\sum_{t = 1}^T \p{u_{A_t}(p_t) -  u_{A_t}(p_t, \, \zeta)} \leq C\zeta^2 \ \text{ for all } \ 0 \leq \zeta < \alpha.
\end{equation}
\end{theo}

Recall that, as the market size gets large, Theorem \ref{thm:main} enables us to estimate gradients of $u_a(p)$
in large-$n$ markets using an amount of randomization that scales as $n^{-\alpha}$ for some $0 < \alpha <  0.5$.
Combined with Theorem \ref{thm:cost}, this result implies that we can in fact estimate gradients of $u_a(p)$
``for free'' via local experimentation when $n$ is large, and that the regret of a platform deploying our platform
matches to first order the regret of an oracle who was able to run first-order optimization on the mean-field limit.

\subsection{Comparison with Rates for Global Experimentation}
\label{sec:bandit}

As discussed above, our local experimentation approach makes two departures from the classical literature
on experimental design under interference, including
\citet{aronow2017estimating}, \citet{athey2018exact}, \citet{baird2018optimal}, \citet{basse2019randomization},
\citet{eckles2017design}, \citet{hudgens2008toward}, \citet{leung2020treatment}, \citet{manski2013identification},
\citet{sobel2006randomized} and \citet{tchetgen2012causal}.
First we use mean-field equilibrium modeling to capture and correct for interference effects;
second, we operationalize our approach in a dynamic setting where a decision maker wants to tune a decision
variable while controlling realized regret while learning.

To highlight the value of mean-field equilibrium modeling, we compare our result from Theorem \ref{thm:learning} to
what can be achieved via the global experimentation baseline that is tailored to sequential decision making, but does
not use equilibrium modeling: Each day $t = 1, \, ..., \, T$, global
experimentation chooses a payment $p_t$ given to all workers on that day, and the observes the corresponding reward $U_t$.
Analogously to the random saturation design discussed in \citet{baird2018optimal} and \citet{hudgens2008toward},
global experimentation does not suffer any bias due to interference because there is no cross-day interference in our model.
The downside of global experimentation is that, unlike our equilibrium modeling based approach, it does not provide the
analyst any direct information about gradients $u'_{A_t}(p_t)$, and this severely limits the ability of global experimentation to effectively
discover a good choice of $p$.

To understand the limits of global experimentation we turn to the literature on continuous-armed bandits
(or zeroth-order optimization), which has established strong lower bounds for closely related problems.
\citet{shamir2013complexity} considers the following setting: We have a sequential decision making problem where,
in each time period, the analyst gets to choose $p_t$ from a bounded interval $I$ and observes a reward $U_t$
with \smash{$\EE{U_t \cond p_t} = u(p_t)$} and \smash{$\Var{U_t \cond p_t} = 1$}; the goal is to choose a sequence $p_t$ that makes the
regret \smash{$\sum_{t = 1}^T (u(p^*) - u(p_t))$} small, where $p^*$ is the maximizer of $u(\cdot)$ over the interval $I$.
\citet{shamir2013complexity} then shows that, even if $u(\cdot)$ is strongly concave, no algorithm can achieve
expected regret that grows slower than $\sqrt{T}$; and, in fact, this result holds even if $u(\cdot)$ is known a priori to be
a quadratic with unit curvature. Further results in this line of work are given in \citet{bubeck2017kernel}. We also note
closely related results by \citet{keskin2014dynamic} who establish a $\sqrt{T}$ lower bound on regret for pricing under
a linear demand model (note that, with linear demand, the seller's profit is quadratic),
and by \citet{nambiar2019dynamic}, who propose a global experimentation scheme driven by random
perturbations to $p_t$ that could be used to achieve $\sqrt{T}$ regret in our model.

The upshot is that, when the daily reward functions $u_{A_t}(p_t)$ are strongly concave and there is meaningful
cross-day noise due to $A_t$, our approach can achieve cumulative regret on the order of $\log(T)$ (corresponding
to a $1/t$ rate of decay in errors), whereas global experimentation cannot improve over $\sqrt{T}$ regret
(corresponding to a $1/\sqrt{t}$ rate of decay in errors).
In other words, our ability to use mean-field modeling to leverage small-scale payment variation within (rather
than across) time periods enables us to fundamentally alter the difficulty of the problem of learning the optimal $p$, and to
improve our rate of convergence in $T$.

Finally, we note that the well-known slow rates of convergence for continuous-armed bandits have led some
authors to studying a query model where we can evaluate the unknown functions $u_{A_t}(\cdot)$ twice
rather than once; for example, \citet{duchi2015optimal} show that two function evaluations can result in substantially faster
rates of convergence than one. The reason for this gain is that, given two function
evaluations, the analyst directly cancel out the main effect of the global noise term $A_t$.
In our setting, it is implausible that a platform
could carry out such paired function evaluations in practice unless, e.g., they simultaneously run experiments
across two identical twin cities. But in this paper, we found that---by leveraging structural information
and mean-field modeling---local experimentation can be used to obtain similar gains over zeroth-order
optimization as one could get via twin evaluation.

\section{Generalizations and Limitations}
\label{sec:generalizations}

So far, we have focused our discussion on a specific a model of a centralized market for freelance labor;
but, as outlined in the introduction, we expect the general principles outlined here to be
more broadly applicable. A full theory of experimental design powered by mean-field equilibria
is beyond the scope of this paper. In this section, however, we take a first step towards a more general
theory by presenting two problem settings of considerable practical interest that are amenable to our approach,
risk-averse suppliers and surge pricing, and discuss another problem, immunization via vaccines, that does
no appear to be amenable to it.

\subsection{Equilibrium Modeling via Generalized Earning Functions}
\label{sec:GEF}

In our motivating model for freelance labor, we considered a setting where
the platform first chooses a distribution $\pi$ and then, for each supplier $i$, draws $P_i \sim \pi$ and promises
to pay the supplier $P_i$ per unit of demand served; the supplier computes $q_A(\pi)$, the expected number
of units of demand they will get to serve if they join the market; finally, each supplier compares their expected revenue
$P_i q_A(\pi)$ to their outside option and chooses whether or not to join the marketplace.
Our main results were that:
1) In large markets, we can unobstrusively estimate a marginal response function via local experimentation;
2) The behavior of this marketplace can be characterized by a mean-field limit;
3) In the mean-field limit, we can transform estimates of the marginal response function into predictions
of the effect of policy-relevant interventions. Thus, in large markets, we can use local experimentation
for optimizing platform choices.

Here, we briefly discuss how to extend our approach to allow for risk-averse suppliers and surge pricing.
In order to do so, we first define choice models for
both problems below, and write down balance conditions generalizing \eqref{eq:choice}. Afterwards,
we conjecture the existence and form of a mean-field equilibrium, and show that the conjectured equilibrium
model lets us again map from consistent estimates of a marginal response function to relevant
counterfactual predictions---using the same recipe as deployed in the rest of this paper.
As discussed further below, what enables us to extend our discussion to these problems is that,
in both cases, we can explain the choices of suppliers in terms of a unifying formalism we refer to
as generalized earning functions.

\begin{exam}[Risk Aversion]
Under risk aversion, supplier utility functions may not scale linearly with their revenue, and
instead there is a concave function $\beta$ such that the relevant quantity for understanding
the suppliers' choices is the expectation of $\beta(\text{revenue})$ \citep{holt2002risk,pratt1978risk}.
Suppose that $\beta(0) = 0$, and that each worker can serve 0 or 1 units of demand.\footnote{Generalizations
to workers who can serve many units of demand are immediate, at the expense of more involved notation.}
Then our balance condition \eqref{eq:choice} becomes
\begin{equation}
\label{eq:risk}
\mu\nth_a(\pi) = \PP[\pi]{Z_i=1 \bbar A=a} = \EE[\pi]{f_{B_i}\p{\beta(P_i) q\nth_a(\mu\nth_a(\pi)))}  \cond A=a }. 
\end{equation}
The curvature of the function $\beta(\cdot)$ thus corresponds to the degree of a supplier's  risk aversion,
and setting $\beta(p)=p$ recovers our original risk-neutral model. 
\end{exam}

\begin{exam}[Supply-Side Surge Pricing]
Several prominent ride sharing platforms deploy surge pricing where, in case of heavy demand, the
platform applies a multiplier (generally greater than $1$) to the original payment in order to encourage
higher supplier participation \citep{cachon2017role,hall2015effects}. As a simple model, suppose
that surge is triggered automatically based on the supply-demand ratio, i.e., there is a
function $s: \RR_+ \rightarrow \RR_+$ such that, in each period, the $i$-th supplier gets paid
$ s(D/T) P_i$ per unit of demand served. Suppliers can anticipate surge and, as in the rest of the
paper, they make decisions based on limiting values of
all random variables. Thus, suppliers anticipate payments $s(d_a / \mu_a(\pi)) P_i$, resulting in a
balance condition
\begin{equation}
\label{eq:surge}
\mu\nth_a(\pi) = \PP[\pi]{Z_i=1 \bbar A=a} = \EE[\pi]{f_{B_i}\p{s\p{\frac{d_a}{\mu\nth_a(\pi)}} P_i \, q\nth_a(\mu\nth_a(\pi)))}  \cond A=a }, 
\end{equation}
where again $s(x) = 1$ recovers our original model.
\end{exam}

In both examples above, we conjecture that---in analogy to Lemma \ref{lem:limit_mup}---a mean-field limit exists and that it can be characterized by analogues
of \eqref{eq:risk} and \eqref{eq:surge} but without the $n$-superscripts. In this case, we can write both mean-field limits in a
unified form via \emph{generalized earning functions}, $\theta: \rp^2 \to \rp$, so that the asymptotic balance condition is
\begin{equation}
\label{eq:GEF}
\mu_a(\pi) = \PP[\pi]{Z_i=1 \bbar A=a} = \EE[\pi]{f_{B_i}\p{\theta(P_i,  \, q_a(\mu_a(\pi)))}  \cond A=a }. 
\end{equation}
In the case of \eqref{eq:risk}, we have $\theta_{risk}(p, \, q) = \beta(p)q$. Meanwhile, for \eqref{eq:surge}, recall that
in the mean-field limit the matching of supply and demand is characterized by the identity
$q_a(\mu_a(\pi)) = \omega(d_a/\mu_a(\pi))$. Thus, our conjecture means that \eqref{eq:surge} converges
to \eqref{eq:GEF} with generalized earning function $\theta_{surge}(p, \, q) = p q s(\omega^{-1}(q))$.

We close this section by carrying out ``step 3'' of the analysis outlined in the first paragraph of this section,
 i.e., by showing
how \eqref{eq:GEF} lets us map from a marginal response function to utility gradients with respect to surge;
we leave verification of the conjectured convergence to \eqref{eq:GEF} for further work.
To this end, fix $a \in \acal$.  First, it is not difficult to show that the changes caused by the introduction of
$\theta(\cdot)$ affects the computation of utility derivative $u_a'(p)$ only through the expression for $\mu_a'(p)$
(cf.~the proof of Proposition \ref{prop:utilityderiv}). Hence, we here only focus on expressions for $\mu_a'(p)$.

Now, we can directly check that a reduced form expression as in \eqref{eq:local_randomization} allows us
to estimate the following marginal response function via local randomization,
\begin{equation}
\label{eq:MRF_GEF}
\Delta_a(p)  = (\nabla\theta)_1(p,q_a(\mu_a(p)))  \EE{f_{B_1}'(\theta(p, q_a(\mu_a(p))) \cond A=a}. 
\end{equation}
where $(\nabla\theta)_i(\cdot, \cdot)$ denotes the $i$-th coordinate of the gradient of $\theta$.
Meanwhile, an argument based on the chain rule similar to that in the proof of Lemma \ref{lem:limit_mup} shows that
\begin{equation}
\mu_a'(p) = {\Delta_a(p) } \, \bigg/  \p{1+ \frac{(\nabla\theta)_2(p, \, q_a(\mu_a(p)))}{(\nabla\theta)_1(p, \, q_a(\mu_a(p)))} \
\omega'\p{\frac{d_a}{\mu_a(p)}} \frac{d_a}{\mu_a^2(p)} \, \Delta_a(p)}. 
\label{eq:GEFdmu}
\end{equation}
Note, furthermore, that all quantities in \eqref{eq:GEFdmu} except $\Delta_a(p)$ are either known a-priori or
can be estimated via observed averages.
The upshot is that the mean-field equilibrium characterized by \eqref{eq:GEF} enables us to map an easy-to-estimate
marginal response function to $\mu_a'(p)$ via \eqref{eq:GEFdmu}. These estimates of $\mu_a'(p)$ can then be directly
used to compute utility gradients $u_a'(p)$ that can be used for first-order optimization.

\subsection{Interference and Choice Modeling}
\label{sec:choice}

Although our approach to interference via equilibrium modeling provides useful insights in many problems of
interest, it does not unlock all problems where we want to understand the effects of deploying an intervention
at scale in a large system. One prominent example to which our approach does not (at least obviously) apply pertains
to the study of vaccine effectiveness in the presence of herd immunity \citep{hudgens2008toward,ogburn2017vaccines}.

\begin{exam}[Vaccine Effectiveness]
\label{exam:vaccine}
We are considering whether to enact a policy that would increase vaccination rates against a contagious disease
in a population where only a moderate fraction of people are currently vaccinated.
Due to the interaction among people within the same geographical vicinity, the risk of infection for any given individual not only
depends on whether they are vaccinated themselves, but also on the overall fraction of infected individuals in the ambient
population (which in turn is modulated by the overall fraction of vaccinated individuals).
Thus simple randomized controlled trials cannot be used to consistently estimate the effect of policies that increase the
overall rate of vaccination, and instead methods that explicitly account for interference are required.
\end{exam}

The classical way to think about experiments for community-level vaccine immunity is to randomize the fraction of people vaccinated
across different disjoint (and thus non-interfering) communities \citep{baird2018optimal,hudgens2008toward}. This approach
is directly analogous to the global experimentation baseline considered throughout this paper, and naturally leads to the question
of whether our approach could be used to design more powerful alternatives.

In analogy to notation used in the rest of the paper, index communities by
$t$ and people within communities by $i$, and let $Z_{it} \in \cb{0, \, 1}$ denote whether the $i$-th person in the $t$-th community
gets infected. We write $\mu(p)$ for the expected fraction of people who get infected in a community in which a fraction $p$ of people
are vaccinated, and focus on estimating $d\mu(p)/dp$, i.e., the decrease in the overall infection rate that can be achieved by increasing
the vaccination probability. In this context, global experimentation seeks to estimate $\mu(p)$ by randomly assigning a single vaccination
probability $p_t \in [0, \, 1]$ to each community, so that each person in community $t$ gets (randomly) vaccinated with probability $p_t$.
In contrast, a local experimentation might consider using individualized randomization probabilities $p_{it} \in [0, \, 1]$ to get a better
handle on $d\mu(p)/dp$.

At first glance, the problem may not appear so different from our leading example.
It seems plausible that the above model sketch could be formalized in a way
that makes it amenable to mean-field asymptotics. Furthermore, in this setting, using symmetric perturbations
$p_{it} = p_t \pm \zeta \varepsilon_{it}$ and regressing $Z_{it}$ on $\varepsilon_{it}$ should recover a well-defined marginal response
function $\Delta(p)$ that, under regularity conditions, corresponds exactly to what's called the (average) direct effect in the statistics literature
\citep{hudgens2008toward,savje2017average}.

At this point, however, we appear to get stuck: Unlike in the main examples considered in this paper, there does not seem to be a
natural way to map from $\Delta(p)$ to our main quantity of interest, namely $d\mu(p)/dp$. In the case of modeling freelance labor,
we assumed that suppliers only care about expected revenue; thus, once $\Delta(p)$ gave us a handle on how they react to
changes in expected revenue due to the platform directly changing $p_{it}$, we were also able to reason about how they might react
to changes in expected revenue due to changes in marketplace conditions that arise from general equilibrium effects.
In the case of vaccine effectiveness, however, there's no a-priori obvious way to connect the direct effect of vaccinating a specific
person to how the same person will react to a change in the overall fraction of the population that's infected. For example, there is presumably
a positive association between how much different people benefit from the vaccine directly and how much they benefit from it via herd
immunity; however, some people may not be responsive to the vaccine and so have zero direct effect, but will still benefit indirectly from the
vaccine via herd immunity. Thus, there appears to be no way to credibly learn about vaccine effectiveness without considering
exogenous variation in the fraction of the population that's infected.

An interesting conceptual distinction between all the positive examples presented in this paper and the above vaccination example
is that, in the former, interference effects are fundamentally due to choices made by participants in the system. For example, in the
case of our model for freelance labor, interference effects arise because suppliers choose not to participate in marketplaces that
are too congested. In contrast, in the vaccination example, getting sick isn't a choice; it's simply a random event whose probability
can be modulated up or down via different vaccination policies and community-level infection levels. The fact that joining a marketplace
is a choice whereas getting sick is not may not matter from the point of view of mean-field asymptotics; however, making assumptions about
how suppliers make choices is what let us credibly connect $\Delta(p)$ with $d\mu(p)/dp$ and proceed with our approach.
The role of choice versus pure chance in understanding best practices for statistical estimation has been the topic of a longstanding
discussion at the intersection of economics and statistics \citep[e.g.,][]{heckman2001micro,imbens2014instrumental,roy1951some};
and, in this context, our result can be seen as one example where simple choice modeling helps motivate a powerful approach to statistical
inference and learning.

\section{Simulation Results}
\label{sec:simu}

We now consider a more comprehensive empirical evaluation of the performance of local versus
global experimentation, building on in the simulation results of Section \ref{sec:motivation},
and compare mean performance of local experimentation
and global experimentation across 1,000 simulation replications. Local experimentation is run for 200 steps,
exactly as described in Section \ref{sec:motivation}, with a random initialization $p_1 \sim \text{Unif}(10, \, 30)$.
Meanwhile, for global experimentation, we consider a collection of strategies that first randomly draw
payments $p_t \sim \text{Unif}(10, \, 30)$ for the first $1 \leq t \leq T$ time periods, fit a spline to
the data (as in the left panel of Figure \ref{fig:baseline}), and then deploy the learned policy for the remaining
$200 - T$ time periods. We consider the choices $T \in \cb{40, \, 60, \, 80, \, \ldots, \, 200}$.
For both methods, we report both in-sample regret, i.e., the mean utility shortfall relative to deploying
the population-optimal $p^*$ for the $T$ learning periods, as well as future expected regret, i.e.,
the expected utility shortfall from deploying the learned policy $\hp$ after the $T$ learning periods.
For local experimentation, we set \smash{$\hp = 2 \sum_{t = 1}^t t\, p_t \, /\, (T(T+1))$} following Corollary \ref{coro:IID},
whereas for global experimentation we set \smash{$\hp$} to be the output of spline optimization discussed above.

\begin{figure}
\begin{center}
\begin{tabular}{cc}
\includegraphics[width=0.48\textwidth]{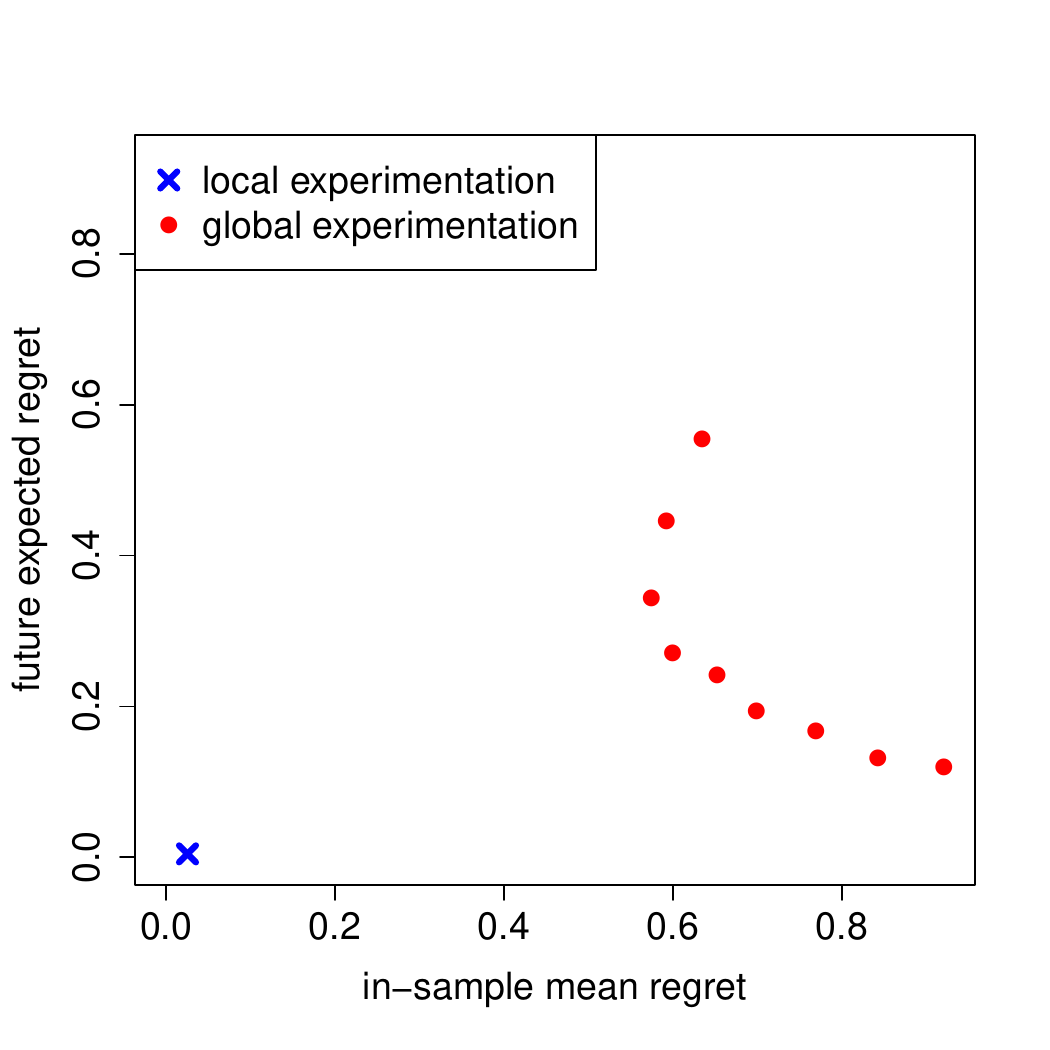} &
\includegraphics[width=0.48\textwidth]{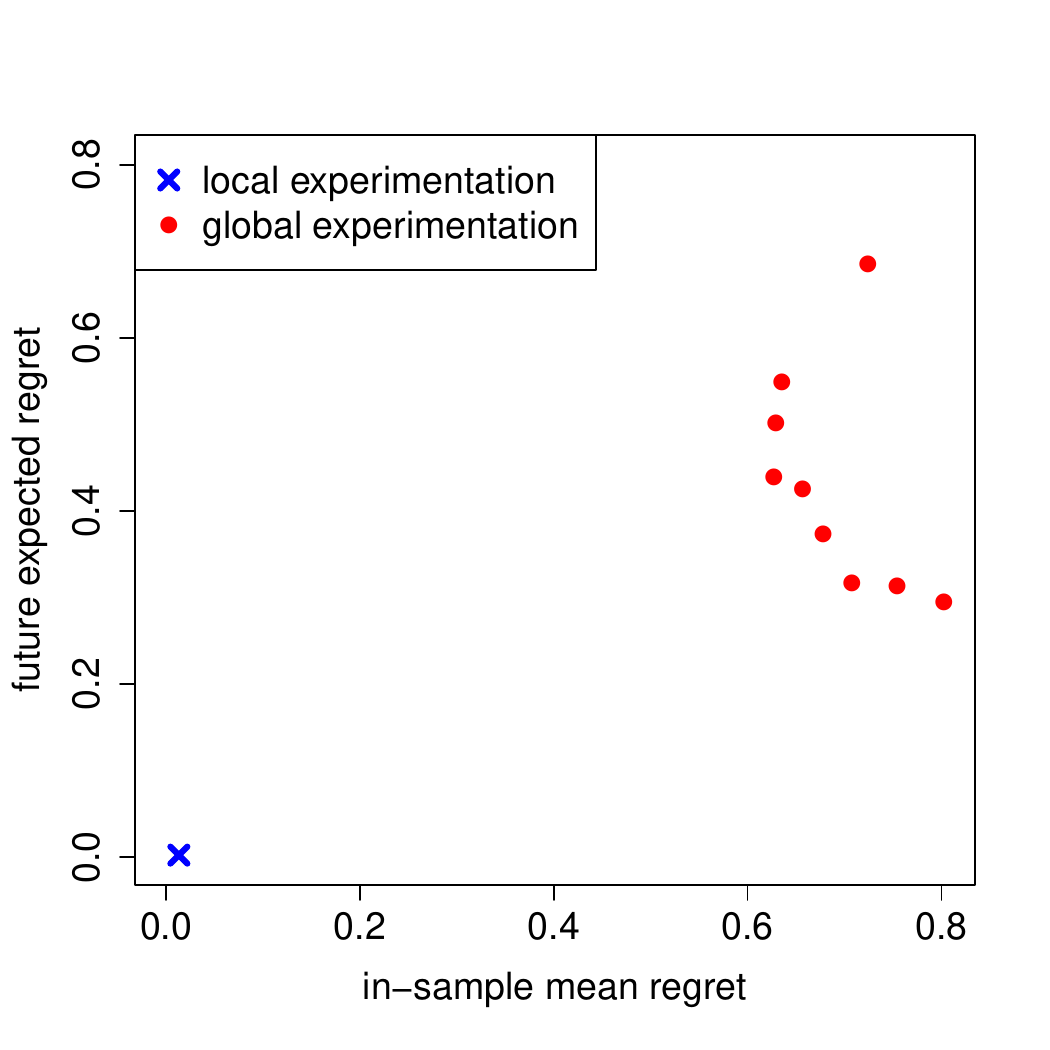} \\
basic setting & supply-side surge
\end{tabular}
\caption{Comparison of the regret of local and global experimentation in the setting of
Section \ref{sec:motivation}, averaged across 1,000 simulation replications. The
global experimentation path is detailed in Section \ref{sec:simu}. In the right panel, the
platform makes a public commitment to multiply supply-side payments by a surge factor \eqref{eq:surge_implementation}.}
\label{fig:simu}
\end{center}
\end{figure}

As seen in the left panel of Figure \ref{fig:simu}, local experimentation outperforms global experimentation
by an order of magnitude along both metrics. Quantitatively, local experimentation achieved mean in-sample
regret of 0.025 and mean future regret of 0.0045. In contrast, the best numbers achieved by global experimentation
for these metrics were 0.57 and 0.12 respectively---and there was not a single choice of tuning parameters
that achieved both. In general, we see that a larger choice of $T$ always improves future regret, whereas for in-sample
regret there is an optimal middle ground that balances exploration and exploitation (here, $T = 80$).

Next, we consider an analogous simulation design, but with supply-side surge pricing. As discussed in
Section \ref{sec:generalizations}, we assume that the platform makes a public commitment to mechanistically
increase supply side payments by a multiplicative factor $s(D/T)$ once the demand $D$ and supply $T$ are realized,
and suppliers take this commitment into account
when choosing whether or not to join the marketplace. Here, we use
\begin{equation}
\label{eq:surge_implementation}
s(D/T) = \frac{D}{T} \, \Big/ \, \omega\p{\frac{D}{T}},
\end{equation}
meaning that, by the properties of $\omega(\cdot)$ as outlined in Definition \ref{def:RAF}, the surge multiplier
is 1 when $D$ is small relative to $T$, but eventually climbs up to the ratio $D/T$ as demand
outpaces supply. This choice of $s(\cdot)$ is by no means optimal; it is simply an example.

As discussed in Section \ref{sec:generalizations}, our analysis of surge relies on a conjecture that relevant
properties of mean-field limits as discussed in Section \ref{sec:limits} still hold with surge.
We work with a limiting platform utility function that depends linearly
on revenue minus costs as in Lemma \ref{lemm:concave_u},
\begin{equation}
u_a(p) = \p{\gamma - p s\p{\frac{d_a}{\mu_a(p)}}} \omega\p{\frac{d_a}{\mu_a(p)}} \mu_a(p).
\end{equation}
As discussed above, we can estimate the $p$-derivative of the expected scaled active supply size, $\mu_a'(p)$, by local
experimentation via \eqref{eq:MRF_GEF} and \eqref{eq:GEFdmu}. Moreover, following the argument of Theorem \ref{thm:main},
we obtain $p$-derivatives of $u_a(p)$ via
\begin{equation}
\begin{split}
&u_a(p) = \p{\gamma - p s\p{\frac{d_a}{\mu_a(p)}}} 
\p{-\omega'\p{\frac{d_a}{\mu_a(p)}} \frac{d_a}{\mu_a(p)} + \omega\p{\frac{d_a}{\mu_a(p)}}} \mu_a'(p) \\
&\ \ \ \ \ \ \ \ - s\p{\frac{d_a}{\mu_a(p)}}\omega\p{\frac{d_a}{\mu_a(p)}}\mu_a(p)
+ p s'\p{\frac{d_a}{\mu_a(p)}} \frac{d_a}{\mu_a(p)} \omega\p{\frac{d_a}{\mu_a(p)}} \mu_a'(p).
\end{split}
\end{equation}
We turn this into a feasible estimator plugging in our local experimentation estimates of $\hmu_a'(p)$
for $\mu_a'(p)$, and estimating the ratio $d_a/\mu_a(p)$ via its sample analogue $D/T$.

Results for learning $p$ are given in the right panel of Figure \ref{fig:simu}. Qualitatively, the results match those
obtained without surge, and local experimentation still outperforms global experimentation by an
order of magnitude. Local experimentation achieved mean in-sample regret of 0.013 and mean future
regret of 0.0024, while the best corresponding numbers achieved by global experimentation
for these metrics were 0.63 and 0.29 respectively.
We also not that adding the automatic surge multiplier as
in \eqref{eq:surge_implementation} decreased the optimal base payment from 17.6 to 15.7,
while increasing optimal mean platform utility by 0.06 (the median utility difference is 0.04).
Thus, in this example, the regret of global experimentation is much larger than the utility gain
from using surge as in \eqref{eq:surge_implementation} relative to not using surge---whereas
the regret of local experimentation is less than the effect of adopting surge.

Finally, we note that the global experimentation baseline considered here---namely our two-phase algorithm that starts
with pure exploration and then moves to pure exploitation---is fairly simple, and it is
possible that a more sophisticated global experimentation baseline could somewhat improve performance.
However one can check that, under reasonable conditions and provided we explore for the first $\sqrt{T}$ periods,
our implemented baseline attains the optimal $\sqrt{T}$ regret rate of \citet{shamir2013complexity} discussed in Section \ref{sec:bandit}.
Thus, more sophisticated methods like Bayesian zeroth-order optimization\footnote{One potentially promising approach
would be to use local experimentation to get gradient estimates $u'_{A_t}(p_t)$, and then incorporate these estimates into
a Bayesian learning framework. It is plausible that this could yield practically meaningfully improvements over the
first-order approach considered in this paper.}
as considered in, e.g., \citet{letham2018constrained} may improve on finite sample performance but cannot improve on the overall regret
rate of our baselines.\footnote{Another class of popular continuous-armed bandit algorithms were introduced by \citet{flaxman2005online}
and \citet{kleinberg2005nearly}. These methods estimate derivatives by noisy function evaluations and then use these
for gradient descent. However, while desirable due to their transparency and ease, these methods suffer cumulative
regret on the order of $T^{3/4}$ in our setting. In our simulations, this class of methods performed worse than the global
experimentation baseline we report results for.}

\section{Discussion}

We introduced a new framework for experimental design in stochastic systems with significant cross-unit interference. The key insight is that, in certain families of models, inference is structured enough to be captured by a small number of key statistics, such as the global demand-supply equilibrium, and the impact of interference can be subsequently accounted for using mean-field and equilibrium modeling. We then proposed an approach based on local experimentation that would allow us to accurately and efficiently estimate the utility gradient in the large-system limit, and use these gradient estimates to perform first-order optimization. 

There are some simplifying assumptions we make in this work that can be relaxed or verified in future research.  For instance, we have assumed that the demand is exogenous. We expect that an extension of our method can be used to capture scenarios where the demand may, for instance, depend on the supply level: For example,  a passenger may be less likely to hail a ride if they know there would be a long wait. Another assumption we made is that the market equilibrium can be reached relatively quickly. While there are recent empirical evidence suggesting that drivers in a ride-sharing platform do respond to payment changes in manner that takes into account the resulting market equilibrium  \citep{hall2019pricing}, it would be interesting to consider a more realistic model where prices may be updated continuously before a new market equilibrium is fully reached. It is less clear how the current model would apply in this setting, which is likely to require a substantially more sophisticated analysis.

We believe that the general approach proposed in this paper, one that leverages mean-field modeling in experimental design, has the potential to be applicable in a wider range of problems.  As one example, we may consider models in which the key statistics that capture the interference patterns are multi-dimensional. This could occur in a marketplace which, instead of being fully centralized, consists of a small number of inter-connected sub-markets.
 For instance, in a ride-sharing platform, the sub-markets may correspond to neighboring cities connected by highways and bridges. 
 In these systems, suppliers' behaviors remain to be primarily  influenced by the local supply-demand equilibrium in  their respective sub-markets. These local equilibria in turn interact with one another due to network effects. Nevertheless, in a large-market regime where the numbers of market participants are relatively large in all sub-markets, while the total number of sub-markets remains the same, 
 we may still use the type of mean-field asymptotics in this paper to account for the interference across both individuals units and sub-markets to efficiently  estimate the effect of payment adjustments. In another direction, we may extend the one-shot equilibrium model adopted in this paper to dynamic settings where the equilibrium emerges gradually gradually according to a stochastic process (e.g., suppliers may adapt to payment variations only over time), and study whether a dynamic version of our mean-field model can be used to analysis the effects of local experimentation in these systems. 
Finally, it would be interesting to investigate whether the local experimentation scheme proposed in this paper can be generalized to estimate higher-order derivatives of the utility function.

\bibliographystyle{plainnat}
\bibliography{references}

\ifx \useplain\undefined
\begin{APPENDICES}
\else
\begin{appendix}
\fi

\normalsize

\section{Proofs of Main Results}
\subsection{Proof of Lemma \ref{lem:sa_unique}}

Fix $a\in \calA$. Recall that all suppliers know the realization of the global state, $a$, and the probability that a given supplier will choose to become active is given by \eqref{eq:choice}. Define 
\begin{equation}
\psi\nth_a(\mu,  \pi) \bydef \EE{ f_{B_1} \lt(P_1 \,  \EE{\Omega\p{D, X}  \bbar A=a}\rt) \cond A=a }, \ \ \ X \sim \text{Binomial}(n, \, \mu).
\end{equation}
That is, $\psi\nth_a(\mu, \pi)$ is the probability of a supplier becoming active under the payment distribution $\pi$, \emph{if} they believe that the active supply size is Binomial $(n, \, \mu)$.  
Note that, since the same payment distribution applies uniformly across all suppliers, so are the probabilities of the suppliers becoming active. As a result, for any $\pi$,  the actual active supply size $T$ will follow a Binomial distribution.  In particular, this implies that $T$ is an equilibrium active supply size if and only if it is binomial with mean $\mu$ that satisfies the following fixed-point equation: 
\begin{equation}
\psi\nth_a(\mu, \pi) = \mu.
\label{eq:balanceq}
\end{equation}
It suffices to show that \eqref{eq:balanceq} admits a unique solution in the domain $\mu \in [0,1]$. Because $f_b(\cdot)$ is by construction non-decreasing, it follows that $\psi\nth_a(\mu, \pi) $ is a continuous function and {non-increasing} in $\mu$: A supplier is  more discouraged from becoming active, if they believe there will be more active suppliers in the market eventually. In particular, the left-hand side of   \eqref{eq:balanceq}, $  \psi\nth_a(\mu, \pi)$, is a non-negative, continuous  and non-increasing function  over $\mu\in [0, 1]$, which implies that \eqref{eq:balanceq} admits a unique solution in $ [0, 1]$. 
\endproof

\subsection{Proof of Lemma \ref{lem:conv_fluid}}

Our argument will leverage the following simple expression for the limiting
derivative of $q(\cdot)$, the proof of which follows immediately from a generalization of a classical result of
\citet{stein1981estimation} to exponential families; the proof is given in Appendix \ref{app:prop:qprimeconvg}. 

\begin{prop}
\label{prop:qprimeconvg}
Fix $a\in \calA$ and $\mu >0$. Then, $\frac{d}{d\mu}q\nth_a(\mu)$ is non-positive, and 
\begin{equation}
\lim_{n \to \infty} \frac{d}{d\mu}q\nth_a(\mu) =\frac{d}{d\mu} \omega\p{{d_a}/{\mu}}= - \omega'(d_a/\mu)\frac{d_a}{\mu^2}. 
\end{equation}
\end{prop}

The claim in \eqref{eq:qconvg} follows directly from the definition of $\Omega$ and Assumption \ref{ass:Omega_tou}, i.e., that $\Omega(d,t) $ converges to $\omega(d/t)$ as $t\to \infty$, and the fact that conditional on $A=a$, $D/n$ concentrates on $d_a$ as $n\to \infty$. 
For \eqref{eq:muconvg}, recall that $\mu\nth_a(p)$ is the solution to the balance equation in \eqref{eq:balanceq}:
$ \mu = \EE{ f_{B_1} (p \, q\nth_a(\mu) ) \cond A=a }$.
By \eqref{eq:qconvg}, and the monotonicity of the functions $f_b(\cdot)$ and $\omega(\cdot)$, we have that $\mu\nth_a(p)$ converges to $\mu_a(p)$ as $n \to \infty$, where $\mu_a(p)$  is the solution to the limiting balance equation given in the statement of Lemma \ref{lem:conv_fluid}. The claim in \eqref{eq:uconvg} follows from \eqref{eq:muconvg}, \eqref{eq:qconvg} and Assumption \ref{ass:Rtou}. 
The convergence of $q^{(n)'}_a(\mu)$ in  \eqref{eq:qprimeconvg} follows from Proposition \ref{prop:qprimeconvg}.
\endproof

\subsection{Proof of Lemma \ref{lem:limit_mup}}

We start by verifying \eqref{eq:mup1}.
By \eqref{eq:balanceq} and the chain rule, we have that 
\begin{align*}
\frac{d}{dp} \mu\nth_a  (p) = & \frac{d}{dp} \psi(\mu\nth_a, \delta_{p}) \\
 =&  \frac{d}{dp}\EE{f_{B_1}\lt(p q\nth_a(\mu\nth_a(p))\rt) \cond A=a } \\
= & \EE{f'_{B_1}\lt( p  q\nth_a  (\mu\nth_a(p)) \rt) \lt( q\nth_a (\mu\nth_a(p)) +  p (q\nth_a)'(\mu\nth_a(p))(\mu\nth_a)'(p) \rt)  \cond A=a}\\
= &  q\nth_a (\mu\nth_a(p)) \EE{ f'_{B_1} ( p  q\nth_a  (\mu\nth_a(p)) )  \cond A=a } \\
&\ \ \ \ \ \ \ \ \ \ \ \  +  p  (q\nth_a)'(\mu_a(p)) (\mu\nth_a)'(p)    \EE{ f'_{B_1} ( p  q\nth_a  (\mu\nth_a(p)) )  \cond A=a} \\
= & \Delta\nth_a(p) +  p   \Delta\nth_a(p)  (q\nth_a)'(\mu\nth_a(p)) (\mu\nth_a)'(p)/q\nth_a(\mu\nth_a(p)).
\end{align*}
The last expression above is linear in $(\mu\nth_a)'(p)$. Re-arranging the equation and solving for $(\mu\nth_a)'(p)$ leads to the desired result. 
Finally, \eqref{eq:deltaconvg} is a direct consequence of Lemma \ref{lem:conv_fluid},
while \eqref{eq:mup2} follows by combining from \eqref{eq:mup1} with \eqref{eq:deltaconvg}
and Lemma \ref{lem:conv_fluid}.
\endproof

\subsection{Proof of Theorem \ref{thm:main}}
\label{sec:pf_thm_main}

The proof will make use of the following two technical results. The first concerns the sensitivity of the system dynamics
with respect to small perturbations $\zeta$, and the second extends calculations from Section \ref{sec:limits} to the utility function
$u\nth_a(p)$. The proofs of these results are given in Appendices \ref{app:prop:perturb_df} and \ref{app:prop:utilityderiv}, respectively. 
Recall that $\mu\nth_a(p, \zeta)$ is the expected fraction of active
suppliers in equilibrium when the payments are $\zeta$-perturbed from $p$, i.e., $\mu\nth_a(p, \zeta) \bydef \EE{T(p, \zeta)/n \bbar A=a}$.  

\begin{prop}
\label{prop:perturb_df}  
Fix $p >0$,  $a\in \calA$ and $n\in \N$. $\mu\nth_a(p, \zeta)$ and $q\nth_a(p,\zeta)$ are twice differentiable functions with respect to $\zeta$, and satisfy: 
\begin{enumerate} 
\item  $ \cb{ \frac{\partial }{\partial \zeta}  \mu\nth_a(p, \zeta) }_{\zeta=0}=  \cb{ \frac{\partial }{\partial \zeta} q\nth_a(\mu\nth_a(p, \zeta))}_{\zeta=0}=0$ for all $n \in \N$. 
\item  There exists $\alpha>0$ such that $ \cb{\frac{\partial^2 }{\partial^2 \zeta}  \mu\nth_a(p, \zeta)}_{\zeta = \zeta_0} $ and $\cb{\frac{\partial^2 }{\partial \zeta^2} q\nth_a(\mu\nth_a(p, \zeta)) } _{\zeta = \zeta_0}$ are bounded uniformly  over all $\zeta_0 \in (0,\alpha) $ and $n \in \N$.   
\end{enumerate}
\end{prop}

\begin{prop}
\label{prop:utilityderiv}
Fix $p>0$ and $a\in \calA$. We have that 
\begin{align}
\label{eq:utility_deriv}
  \limn \frac{d}{d p}  u_a{\nth} (p)   = &  u'_a(p)  \nln
  = &  \mu_a'(p)\left[ r\p{\frac{d_a}{\mu_a(p)}} - p \omega\p{\frac{d_a}{\mu_a(p)}} - \p{ r'\p{\frac{d_a}{\mu_a(p)}} - p \omega'\p{\frac{d_a}{\mu_a(p)}}}\frac{d_a}{\mu_a(p)}\right] \nln
& \ \   - \omega\p{\frac{d_a}{\mu_a(p)}}\mu_a(p), \nnb
\end{align}
where $u_a(\cdot)$ is defined in \eqref{eq:uconvg}.
\end{prop}

Now, recall that we are considering the case where the platform employs an $\eta$-perturbed payment distribution $\pi_{p, \eta}$,
with $P_i = p+\eta \varepsilon_i$ (\eqref{eq:zEta}).  Define the estimators
\begin{align}
\bD =& D/n, \quad \bZ = T/n= \frac{1}{n}\sum_{i=1}^n Z_i , \\
\hDelta = & \zeta_n^{-1} \ {\hCov{Z_i, \, \varepsilon_i}} \,\big/\, {\hVar{\varepsilon_i}}, 
\end{align}
so that $\bD$ and $\bZ$ correspond to the scaled demand and active suppliers, respectively, and $\hDelta$ is the scaled regression coefficient of $Z_i$ on $\varepsilon_i$. Finally, define the estimator
\begin{equation}
 \hUps =  \hDelta/\lt(1+ \frac{p \bD \hDelta \omega'\p{\bD/ \bZ}}{\bZ^2 \omega\p{\bD / \bZ}}\rt). 
\end{equation} 
Our main remaining task is to show that, under the stated conditions,
\begin{align}
\bD \to & d_a, \, \bZ \to  \mu_a(p), \, \hDelta  \to  \Delta_a(p), \, \mbox{and } \hUps    \to  \mu_a'(p), 
\label{eq:estm_convg}
\end{align}
in $L_2$ as $n\to \infty$, for any $a \in \acal$. In light of Proposition \ref{prop:utilityderiv}, the  desired conclusion \eqref{eq:Gamma} then follows immediately by
combining \eqref{eq:Gamma_def} with \eqref{eq:estm_convg} with \eqref{eq:utility_deriv} and invoking Slutsky's lemma.

We now turn to proving \eqref{eq:estm_convg}.
First, we note that the fact that $\hUps \to  \mu_a'(p)$ follows directly by combining the first three convergence claims in \eqref{eq:estm_convg} with Lemma \ref{lem:limit_mup}. The fact that $\bZ \to d_a$ follows from our definition \eqref{eq:Dconv0}. For $\bZ \to \mu_a(p)$, note that by Chernoff bound we know that $\bZ$ concentrates on $\mu\nth_a(p,\zeta_n)$  as $n\to \infty$. Furthermore, we have that
\begin{align}
\limn  \mu\nth_a(p,\zeta_n) = & \limn \mu\nth_a(p) + \limn \p{ \mu\nth_a(p,\zeta_n) -  \mu\nth_a(p) } \nln
\sk{a}{=}& \mu_a(p) +  \limn \p{ \mu\nth_a(p,\zeta_n) -  \mu\nth_a(p) }  \nln
\sk{b}{=} & \mu_a(p), 
\end{align}
where steps $(a)$ and $(b)$ follow from Lemma \ref{lem:conv_fluid} and Proposition \ref{prop:perturb_df}, respectively. Together, this shows that $\bZ \to (p)$ in $L_2$. 

Finally, it remains to show that $\hDelta  \to  \Delta_a(p)$. To this end, we first observe the following fact: There exists a constant $C > 0$ such that, for every $\zeta$ and $n$ 
\begin{equation}
\label{eq:identification}
\abs{ \frac{1}{\zeta} \Cov[n]{Z_i, \, \varepsilon_i \cond A = a} - \Delta\nth_a(p)} \leq \zeta C. 
\end{equation}
To prove \eqref{eq:identification}, note that given $\zeta$-perturbed payments $P_i = p + \zeta \, \varepsilon_i$, we have
$$ \Cov[n]{Z_i, \, \varepsilon_i \cond A = a} = \EE[n]{\varepsilon_i f_{B_i}\p{\p{p + \zeta \varepsilon_i} q_a\nth\p{\mu\nth_a\p{p, \, \zeta}})} \cond A=a }. $$
We can then take the limit $\zeta \rightarrow 0$, and verify that there exists $C>0$ such that
\begin{align*}
& \abs{\EE[n]{\frac{1}{\zeta}\varepsilon_i f_{B_i}\p{\p{p + \zeta \varepsilon_i} q_a\nth\p{\mu\nth_a\p{p, \, \zeta}})} \cond A=a} \right.  \nln
& \ \ \ \ \ \ \  \left. - q_a\nth(p,0) \EE[n]{f_{B_i}'(p q_a\nth(\mu\nth_a\p{p,0})) \cond A=a} } \leq \zeta C 
\end{align*}
Here, we used Proposition \ref{prop:perturb_df}, and specifically the fact that  $\cb{ \frac{\partial}{\partial \zeta} q_a\nth\p{\mu\nth_a\p{p, \, \zeta}}}_{\zeta = 0} = 0$, both $f_{B_i}(\cdot)$ and $q_a\nth(\mu\nth_a\p{p, \, \cdot})$ are  twice differentiable with  bounded second derivatives  uniformly over $n$, and $\varepsilon_i$
has variance 1.  Finally, by Slutsky's lemma and conditionally on $A = a$, we have
$$ \hDelta - \frac{\Cov[n]{Z_i, \, \zeta_n\varepsilon_i \cond A = a}}{\Var[n]{\zeta_n\varepsilon_i \cond A = a}} \rightarrow_p 0. $$
As $\Var[n]{\varepsilon_i \cond A = a} = 1$, we conclude that $\hDelta  \to_p  \Delta_a(p)$
using \eqref{eq:identification} and Lemma \ref{lem:limit_mup}.
\endproof

\subsection{Proof of Theorem \ref{thm:learning}}
\label{sec:pf_learning}

Given the form of \eqref{eq:OMD}, we can use Lemma 1
of \citet{orabona2015generalized} to check that
\begin{equation}
\label{eq:deterministic_bound}
\sum_{t = 1}^T t(p - p_t) \hGamma_t \leq \frac{1}{2\eta} \sum_{t = 1}^T t (p - p_{t})^2 + \frac{\eta}{2} \sum_{t = 1}^T \hGamma_t^2.
\end{equation}
We then can replace the gradient estimates $\hGamma_t$ with their mean-field limits $u'_{A_t}(p_t)$
provided we add appropriate error terms as follows,
\begin{align*}
&\sum_{t = 1}^T  t(p - p_t) u'_{A_t}(p_t)
\leq \frac{1}{2\eta} \sum_{t = 1}^T t(p - p_{t})^2 + \frac{\eta}{2} \sum_{t = 1}^T u'_{A_t}(p_t)^2 \\
&\ \ \ \ +  \sum_{t = 1}^T  t (p - p_t) \p{u'_{A_t}(p_t) - \hGamma_t} +  \frac{\eta}{2} \sum_{t = 1}^T\p{\hGamma_t^2 - u'_{A_t}(p_t)^2}.
\end{align*}
Then, given the result in Theorem \ref{thm:main} we see that, for any $\varepsilon > 0$,
\begin{equation}
\label{eq:gradient_conv0}
\sum_{t = 1}^T  t (p - p_t) u'_{A_t}(p_t)\leq  \frac{1}{2\eta} \sum_{t = 1}^T t (p - p_{t})^2 + \frac{\eta}{2} \sum_{t = 1}^T u'_{A_t}(p_t)^2 + \varepsilon
\end{equation} 
with probability tending to 1 as $n$ gets large. Noting that $\abs{u'_{A_t}(p_t)} < M$, this simplifies to
\begin{equation}
\label{eq:gradient_conv}
\sum_{t = 1}^T  t (p - p_t) u'_{A_t}(p_t)\leq  \frac{1}{2\eta} \sum_{t = 1}^T t (p - p_{t})^2 + \frac{\eta M^2 T}{2}
\end{equation} 
with probability tending to 1.
The desired statement \eqref{eq:consistency} follows by leveraging the remaining assumptions from
the theorem statement: $\sigma$-strong concavity of $u_{A_t}(\cdot)$ implies that
\begin{equation}
\label{eq:strong_conc}
u_{A_t}(p) \leq u_{A_t}(p_t) + (p - p_t) u'_{A_t}(p_t) - \frac{\sigma}{2} (p - p_t)^2,
\end{equation}
and we use the above to replace the left-hand side expression of \eqref{eq:gradient_conv} while noting that $\sigma > \eta^{-1}$.
\endproof

\subsection{Proof of Corollary \ref{coro:IID}}

Let $p^*$ be the maximizer of $u(\cdot)$ over $I = [c_-, \, c+]$.
By \eqref{eq:consistency} we have
\begin{equation}
\begin{split}
&\limn \PP{\frac{1}{T} \sum_{t = 1}^T t \p{u(p^*) -  u(p_t)}  \leq \frac{Z_T}{T} + \frac{\eta M^2}{2}} = 1, \\
&Z_T = \sum_{t = 1}^T t \p{u(p^*) - u(p_t) + u_{A_t}(p^*) -  u_{A_t}(p_t)}.
\end{split}
\end{equation}
Paired with strong concavity of $u(p)$ around $p^*$ and the fact that $u'(p^*) = 0$, this implies
\begin{equation}
\label{eq:consistency_l2_start}
\limn \PP{\frac{\sigma}{2} \frac{1}{T} \sum_{t = 1}^T  t(p^* - p_t)^2
\leq \frac{Z_T}{T} + \frac{\eta M^2}{2}} = 1.
\end{equation}
In order to verify the desired result, our next step is to bound $Z_T$.
First, because $p_t$ is chosen before we get to learn about $A_t$, 
$Z_t$ is a martingale. Second, because the derivative of $u_a(p)$ is uniformly
bounded by $M$, we have $\abs{Z_{t} - Z_{t-1}} \leq 2 M t \abs{p_t - p^*}$ for
all $t$. Thus, using Hoeffding's lemma to bound the moment-generating function
of a bounded random variable, these two facts together imply that
$$ \EE{\exp\p{c (Z_t - Z_{t-1})} \cond p_t, \, \cb{Z_s, \, p_s}_{s = 1}^{t-1}} \leq \exp\p{\frac{1}{2} c^2 M^2  t^2 \p{p_t - p^*}^2}, $$
and so
$$ Y_t = \exp\p{c Z_t - \frac{1}{2} c^2 M^2 \sum_{s = 1}^t s^2 \p{p_s - p^*}^2} $$
is a super-martingale for any $c > 0$. Thus, by Markov's inequality,
\begin{equation}
\label{eq:markov_bound}
\PP{Z_T \geq  \frac{c M^2}{2} \sum_{t = 1}^t t^2 \p{p_t - p^*}^2 + \frac{-\log(\delta)}{c}} \leq \delta.
\end{equation}
for any $0 < \delta < 1$. Pairing \eqref{eq:consistency_l2_start} and \eqref{eq:markov_bound}
with $c = \sigma / (2M^2 T)$ then yields (recall that $\eta > \sigma^{-1}$)
\begin{equation}
\limsup_{n \rightarrow \infty} \PP{\frac{\sigma}{4} \frac{1}{T} \sum_{t = 1}^T  t(p^* - p_t)^2
\leq  \eta M^2 \p{2 \log\p{\delta^{-1}} + 1/2}} \geq 1 - \delta.
\end{equation}
Finally, the desired result follows by noting that
$$ \frac{T^2}{2} \p{p^* - \bar{p}_T}^2 \leq \sum_{t = 1}^T  t(p^* - \bar{p}_T)^2 \leq  \sum_{t = 1}^T  t(p^* - p_t)^2. $$
\endproof

\subsection{Proof of Theorem \ref{thm:cost}}

Using Proposition \ref{prop:perturb_df} and a first-order Taylor expansion with a Lagrange-form remainder,
we immediately see that there is a $C > 0$ such that, for all $n \geq n_0$ and $0 \leq \zeta < \alpha$,
$$ \sum_{t = 1}^T \p{u\nth_{A_t}(p_t) -  u\nth_{A_t}(p_t, \, \zeta)} \leq CT\zeta^2. $$
Since this bound holds for all $n \geq n_0$, it also holds in the limit $n \rightarrow \infty$, thus implying \eqref{eq:cost}.

\section{Proofs of Technical Results}

\subsection{Proof of Lemma \ref{lemm:concave_u}}
\label{sec:concave_u}

Fix $a\in \calA$. It follows from Lemma \ref{lem:conv_fluid} that 
\begin{equation}
u_a (p) = \gamma \mu_a(p) q_a(\mu_a(p)) - p \mu_a(p) \omega(d_a/p)  = (\gamma-p) \mu_a(p)q_a(\mu_a(p)), 
\end{equation}
where 
\begin{equation}
q_a(\mu) = \omega(d_a/\mu). 
\end{equation}
We have that
\begin{align}
u_a''(p) = (\gamma - p) \frac{d^2}{dp^2}\p{q_a(\mu_a(p)) \mu_a(p) }- 2 \frac{d}{dp}\p{q_a(\mu_a(p)) \mu_a(p)}. 
\end{align}
Note that $q_a(\mu_a(p)) \mu_a(p) $ is the normalized amount of demand that ends up being served. The next result shows that $q_a(\mu_a(p)) \mu_a(p) $ is non-decreasing in the payment $p$; The proof is given in Appendix \ref{app:prop:demincreasp}. 
\begin{prop} 
\label{prop:demincreasp}
There exists $c>0$, such that
\begin{equation}
\frac{d}{dp}\p{q_a(\mu_a(p))\mu_a(p)} \geq c, \quad \mbox{for all $p\in (c_0, \gamma)$. }
\end{equation}
\end{prop}

Because $p<\gamma$, in light of Proposition \ref{prop:demincreasp}, in order to show that  $u_a(\cdot)$ is strictly concave, it suffices to demonstrate that 
\begin{equation}
 \frac{d^2}{dp^2}\p{q_a(\mu_a(p)) \mu_a(p) } \leq 0. 
\end{equation}
To this end, we have that
\begin{align}
&  \frac{d^2}{dp^2}\p{q_a(\mu_a(p)) \mu_a(p) }  \nln
 =&  \mu_a''(p) \p{\mu_a(p) q'\p{\mu_a(p)}+ q_a(\mu_a(p))} + \mu_a'(p)^2 \p{\mu_a(p) q''(\mu_a(p)) + 2q_a'(\mu_a(p))}\nln
\sk{a}{=} & \mu_a''(p) \p{\mu_a(p) q'\p{\mu_a(p)}+ q_a(\mu_a(p))}+
  \mu_a'(p)^2 \frac{d_a^2 \omega''(d_a/\mu_a(p))}{ \mu_a(p)^3}  \nln
\sk{b}{=}& \mu_a''(p) \p{ \omega(d_a/\mu_a(p)) - \omega'(d_a/\mu_a(p))\frac{d_a}{\mu_a(p)}}+
  \mu_a'(p)^2 \frac{d_a^2 \omega''(d_a/\mu_a(p))}{ \mu_a(p)^3} 
  \label{eq:d2qmmu}
\end{align}
where steps $(a)$ and $(b)$ follow from the fact that $q_a(\mu) = \omega(d_a/\mu)$. Because $\omega(\cdot)$ is concave, it follows that the second term in \eqref{eq:d2qmmu} is non-positive. Furthermore, recall from Definition \ref{def:RAF} that $\omega(\cdot)$ is  concave and $\omega(0)=0$. By Proposition \ref{prop:weakconcave}, we have that 
\begin{equation}
 \omega(d_a/\mu_a(p)) - \omega'(d_a/\mu_a(p))\frac{d_a}{\mu_a(p)} \geq 0. 
 \end{equation} 
In the remainder of the proof, we will focus on showing that 
\begin{equation}
\mu_a''(p)  \leq 0, 
\end{equation}
which would imply the strong concavity of $u_a(\cdot)$. 

Recall that, by construction, the average choice function $f_a(\cdot)$ is non-decreasing and concave.
Recall from \eqref{eq:balanceq} that $\mu_a(p)$ satisfies the fixed-point equation: 
\begin{equation}
\mu_a(p) = \psi_a(\mu_a(p), p) 
\label{eq:mubalan2}
\end{equation}
where
\begin{equation}
\psi_a(\mu, p) = f_a(p q_a(\mu)). 
\end{equation}
Twice-differentiating \eqref{eq:mubalan2} with respect to $p$, we obtain that 
\begin{align}
\mu_a''(p) = \cb{\frac{\partial }{\partial \mu }\psi_a(\mu, p)}_{\mu = \mu_a(p)} \mu_a''(p) + (\mu_a'(p), 1)  \bfH_{\psi_a}(\mu_a(p),p)(\mu_a'(p), 1)^\intercal, 
\end{align}
where $\bfH_{\psi_a}(\cdot, \cdot)$ denotes the Hessian of $\psi_a(\cdot)$. This leads to 
\begin{equation}
\mu_a''(p)  = \p{1-  \cb{\frac{\partial }{\partial \mu }\psi_a(\mu, p)}_{\mu = \mu_a(p)}}^{-1}  (\mu_a'(p), 1)   \bfH_{\psi_a}(\mu_a(p),p) (\mu_a'(p), 1) ^\intercal. 
\end{equation}
Note that since $f$ and $q$ are non-increasing, we have that $ \cb{\frac{\partial }{\partial \mu }\psi_a(\mu, p)}_{\mu = \mu_a(p)} \leq 0$.  It remains to verify that 
\begin{equation}
(\mu_a'(p), 1)   \bfH_{\psi_a}(\mu_a(p),p) (\mu_a'(p), 1) ^\intercal \leq 0. 
\label{eq:HinnerNeg}
\end{equation}
Recall that $\psi_a(\mu, p) = f_a(p q_a(\mu))$. We have that
\begin{align}
\frac{\partial^2}{\partial \mu^2} \psi_a(\mu, p) =&  f''_a(pq_a(\mu)) (pq_a'(\mu))^2 + f'_a(pq_a(\mu))pq_a''(\mu),  \nln
\frac{\partial^2}{\partial p^2} \psi_a(\mu, p) =&  f''_a(pq_a(\mu))q_a(\mu)^2 \nln
\frac{\partial^2}{\partial p\partial \mu} \psi_a(\mu, p) =&  f''_a(pq_a(\mu)) pq_a'(\mu) q_a(\mu) + f'_a(pq_a(\mu)) q_a'(\mu). 
\end{align}
Rearranging terms, and using the fact that $q_a(\mu) = \omega(d_a/\mu) $,  we can decompose $\bfH_{\psi_a}(\mu,p) $ as follows: 
\begin{equation}
\bfH_{\psi_a}(\mu,p) = f''_a(p q_a(\mu)) \bfA  + \frac{ f'_a(p q_a(\mu) ) \omega''(d_a/\mu) d_a^2}{\mu^4} \bfB +\frac{ f'_a(p q_a(\mu)) \omega'(d_a/\mu)d_a} {\mu^2}\bfC
 \label{eq:Hpsi}
 \end{equation} 
 where 
 \begin{align}
 \bfA = & (pq_a'(\mu), \, q_a(\mu)) \otimes(pq_a'(\mu) , \,  q_a(\mu)), \nln
  \bfB = & \begin{bmatrix}
1 & 0 \\
0 & 0
\end{bmatrix},  
\quad 
  \bfC =  \begin{bmatrix}
\frac{2p}{\mu} & -1 \\
-1 & 0
\end{bmatrix}. \nnb
 \end{align}
We make the following observations concerning the three terms in \eqref{eq:Hpsi}. For the first term, $\bfA$ is the outer product of $ (pq_a'(\mu), \, q_a(\mu))$ with itself and is hence positive semi-definite. Since $f_a(\cdot)$ is concave and  hence $f''_a(\cdot)<0$, we have that $ f''_a(pq_a(\mu))  \bfA$ is negative semi-definite, i.e., 
\begin{equation}
f''_a(pq_a(\mu))  \bfA \preceq 0. 
\label{eq:Anegsem}
\end{equation}
For the second term, note that $f'_a(pq_a(\mu)) >0$ and $\omega''(\cdot)<0$ due to the concavity of $\omega(\cdot)$. Therefore, we have that 
\begin{equation}
 \frac{ f'_a(p q_a(\mu) ) \omega''(d_a/\mu) d_a^2}{\mu^4} \bfB  \preceq 0. 
 \label{eq:Bnegsem}
\end{equation}
For the third term, since we are only interested in the properties of $\bfC$ along the specific direction $(\mu_a'(p), \, 1)$, it suffices to show that when $\mu = \mu_a(p)$, $(\mu_a'(p), 1)   \bfC (\mu_a'(p), 1) ^\intercal $ is non-positive. This claim is isolated in the form of the following proposition; The proof is given in Appendix \ref{app:prop:Cmu'1neg}. 
\begin{prop}
\label{prop:Cmu'1neg}
\begin{equation}
(\mu_a'(p), 1)   \bfC (\mu_a'(p), 1) ^\intercal \leq 0. 
\label{eq:Cmu'1neg}
\end{equation}
\end{prop}
By combining \eqref{eq:Anegsem}, \eqref{eq:Bnegsem} and Proposition \ref{prop:Cmu'1neg}, we have proven  \eqref{eq:HinnerNeg}, i.e., 
\begin{equation}
(\mu_a'(p), 1)   \bfH_{\psi_a}(\mu_a(p),p) (\mu_a'(p), 1) ^\intercal \leq 0. 
\end{equation}
This in turn proves Lemma \ref{lemm:concave_u}.   \qed

\subsection{Proof of Proposition \ref{prop:weakconcave} }
\label{app:prop:weakconcave}
 We have that 
\begin{align}
g(x)  = & g(0) + \int_{0}^x g'(s) ds \sk{a}{\geq}  \int_{0}^x g'(x) ds  \sk{b}{\geq}  x g'(x), 
\end{align}
where $(a)$ follows from the assumption that $g(0)\geq0$, and $(b)$ from the concavity of $g(\cdot)$, which implies that $g'(\cdot)$ is non-increasing  over $x >0$.
For the second statement, we first note that if, for some $c > 0$,
\begin{equation}
\label{eq:check}
g'(x) \geq c \ \frac{2}{x} g\p{\frac{x}{2}},
\end{equation}
then we have that
$$ g(x) = g\p{\frac{x}{2}} + \int_{x/2}^x g'(s) \ ds \geq  g\p{\frac{x}{2}} + \frac{x}{2} g'(x) \geq \p{1 + c} g\p{\frac{x}{2}}. $$
We then conclude by arguing by contraction. Suppose that, for each $M \geq 0$, there exists some $x \geq M$
satisfying \eqref{eq:check}; then, by the above argument and noting that $g(x)$ is non-negative and concave (and thus non-decreasing),
we must either have $g(x) = 0$ for all $x$, or $\lim_{x \rightarrow \infty} g(x) = \infty$. Thus, under our stated
assumptions, the condition \eqref{eq:check} can only hold on a finite interval for any value of $c > 0$, and so
$$ \limsup_{x \rightarrow \infty} \frac{x g'(x)}{g(x)} \leq \limsup_{x \rightarrow \infty}  2 \frac{(x/2) g'(x)}{g(x/2)} = 0. $$
\qed

\subsection{Proof of Proposition \ref{prop:qprimeconvg}}
\label{app:prop:qprimeconvg}

We start by verifying a useful property that applies to any exponential family with discrete support.

\begin{defi}
\label{def:expofam}
Let $\{X\}$ a family of discrete random variables and parameterized by $\theta \in \Theta \subset \R$. We say that $\{X\}$ is an exponential family, if the probability mass function (PMF) $f_\theta$ for $X$ can be expressed as
\begin{equation}
f_\theta(x) = h(x)\exp(\eta(\theta) T(x) - A(\theta)), \quad x\in \mathbb{Z}.
\end{equation}
where $T(\cdot)$ is referred to as the sufficient statistic, and $\eta(\theta)$ the natural parameter. 
\end{defi}

We have the following identity, which is a simple generalization of a result proved by \citet{stein1981estimation}
for Gaussian random variables.

\begin{lemm}
\label{lem:SteinIden}
Fix an exponential family of random variables $X$ with discrete support $\xx$ parametrized by $\theta\in \Theta$ defined over a finite subset of $\R$, with sufficient statistic $T(X)$ and natural parameter $\eta(\theta)$. Then,  for any function $g: \xx \to \R$, we have that
\begin{equation}
	\frac{d}{d\theta}\EE[\theta]{g(X)} = \Cov[\theta]{g(X), \, T(X)}, 
\end{equation}
for all $\theta$ in the interior of $\Theta$. 
\end{lemm}

\bpf  First, observe that $\sum_{x} f_\theta(x) = h(x)\exp(\eta(\theta) T(x) - A(\theta)) = 1$. Taking derivatives on both sides with respect to $\theta$, we obtain that 
\begin{align}
0 = & \sum_{x} h(x)\p{T(x)\eta'(\theta) - A'(\theta)}\p{\exp(\eta(\theta)T(x) - A(\theta)}  \nln
=& \EE[\theta]{T(X)}\eta'(\theta) - A'(\theta)
\end{align}
which implies that 
\begin{equation}
\EE[\theta]{T(X)} =  A'(\theta)/\eta'(\theta). 
\label{eq:ET(Xthe)}
\end{equation}
We have that
\begin{align}
	\frac{d}{d\theta}\EE[\theta]{g(X)} =& \sum_{x \in \xx} g(x) \frac{\partial}{\partial \theta} f_\theta(x) \nln 
	=& \sum_{x \in \xx} g(x) \p{T(x)\eta'(\theta) - A'(\theta)}h(x) \exp\p{\eta(\theta)T(x) - A(\theta)} \nln
	=& \eta'(\theta) \EE[\theta]{g(X) \p{T(X) - A'(\theta)/\eta'(\theta)}} \nln
	\sk{a}{=}& \eta'(\theta) \EE[\theta]{g(X) \p{T(X) - \EE[\theta]{T(X)}}} \nln
	\sk{b}{=}& \eta'(\theta) \EE[\theta]{\p{g(X)-\EE{g(X)}} \p{T(X) - \EE{T(X)}}} \nln
	=&  \Cov[\theta]{g(X), \, T(X)}, 
\end{align}
where $(a)$ follows from \eqref{eq:ET(Xthe)}, and $(b)$ from the fact that $\EE[\theta]{ T(X) - \EE[\theta]{T(X)}}=0$.  This proves Lemma \ref{lem:SteinIden}. \qed

We are now ready to prove the stated result. We will assume that all probabilities are calculated by conditioning on $A=a$, and thus omit it from our notation. 
The fact that $\frac{d}{d\mu}q\nth_a(\mu)$ is non-positive follows directly from the fact that $\Omega(d,t)$ is non-increasing in $t$ (Assumption \ref{ass:Omega_tou}). 
The PMF of an $(n, \, \mu)$ Binomial random variable $X$ can be written as 
\begin{equation}
f_\mu(x) = {n \choose x} \exp\p{x \log\lt(\frac{\mu}{1-\mu}\rt) + n\log(\mu(1-\mu))}. 
\end{equation} 
In particular, the set of Binomial random variables forms an exponential family, with natural parameter $\eta(\mu) = \log\lt(\frac{\mu}{1-\mu}\rt)$ and sufficient statistic $T(X) = X$.  
We now employ Lemma \ref{lem:SteinIden} above. Define 
\begin{equation}
H \bydef X - \EE[\mu]{X} = X - n\mu. 
\end{equation}
For a fixed $d$, we have that{\footnote{The notation $x \in a \pm b$ denotes $x \in [a-b, a+b]$. }}
\begin{align}
\frac{d}{d\mu} \EE[\mu]{\Omega(d, X)} (\eta'(\mu))^{-1}
\sk{a}{=} & \Cov[\mu]{\Omega(d, X), \, X} \nln
 = & \EE[\mu]{\Omega(d, X)H } \nln
 \sk{b}{=}  & \EE[\mu]{\omega(d/X)H} + \EE[\mu]{l(d,X)H}\nln
\sk{c}{\in} & \EE[\mu]{\omega(d/X)H} \pm \sqrt{\EE[\mu]{l(d,X)^2}\EE{H^2}}\nln
=  & \Cov[\mu]{\omega(d/X), X} \pm \sqrt{\EE[\mu]{l(d,X)^2}\EE{H^2}}
\end{align}
where $(a)$ follows from Lemma \ref{lem:SteinIden}, $(b)$ from Assumption \ref{ass:Omega_tou}, and $(c)$ from the Cauchy–Schwarz inequality. Taking expectation with respect to $d\sim D$ on both sides, we obtain 
\begin{equation}
\frac{d}{d\mu} \EE[\mu]{\Omega(D, X)} (\eta'(\mu))^{-1} \in  \EE[\mu]{\omega(D/X)} \pm \sqrt{\EE[\mu]{l(D,X)^2}\EE{H^2}}. 
\label{eq:dmuQdX}
\end{equation}
We next bound each of the two terms in \eqref{eq:dmuQdX}. For the second term, recall that $l(\cdot)$ is bounded and $|l(d,t)| = o(1/\sqrt{d}+1/\sqrt{t})$ (Assumption \ref{ass:Omega_tou}). Furthermore, it follows from the Chernoff bound and \eqref{eq:Dconv}, respectively, that 
\begin{equation}
\PP[\mu]{X \geq n\mu/2}, \ \  \PP{D \geq n d_a/2} =o(1/n) \mbox{ as $n\to \infty$}. 
\end{equation}
This implies that 
\begin{equation}
\EE[\mu]{l(D,X)^2} = o(1/n),  \mbox{ as $n\to \infty$}. 
\end{equation}
Furthermore, note that
\begin{equation}
 \EE[\mu]{H^2} = \Var[\mu]{X} = \mcal{O}(n).
 \end{equation} 
Combining the above two equations, we conclude that 
\begin{equation}
 \sqrt{\EE[\mu]{l(D, X)^2}\EE{H^2}}  = o(1), \quad \mbox{as $n \to \infty$}. 
 \label{eq:ElDX}
\end{equation}
Next, we turn to the first term in \eqref{eq:dmuQdX}, which will follow from the following result.  Fix $\delta \in (0,\mu)$, and define the event 
\begin{equation}
\calE = \{ |H/n| < \delta \}. 
\end{equation}
Using Taylor expansion  on the function  
\begin{equation}
h(x) \bydef \omega\p{\frac{d/n}{\mu + x}} 
\label{eq:hfctdef}
\end{equation}
and the smoothness of $\omega$, we have that there exists a constant $c_1>0$ such that, for all $n$ and $d$, 
\begin{align}
h(x) \in   h(0) + h'(0)x \pm c_1x^2, \quad \forall x \in  [-\delta, \delta]. 
\label{eq:hfcttaylor}
\end{align}
Fix $d \in \rp$.  We have that 
\begin{align}
  \EE[\mu]{\omega\p{{d}/{X}}H }   = &  \EE[\mu]{  \omega\p{\frac{d/n}{(n \mu  + H)/n}}H } \nln
 = &  \EE[\mu]{  h(H/n) H } \nln
 = &  \EE[\mu]{  \iden{\calE}h(H/n) H } +  \EE[\mu]{  \iden{\overline{\calE}}h(H/n) H } \nln
  \sk{a}{\in} &  \EE[\mu]{  \iden{\calE}h(H/n) H } \pm c_2\sqrt{\pb(\overline{\calE}) \EE[\mu]{H^2}} \nln
    \sk{b}{\in} &  \EE[\mu]{  \iden{\calE}h(H/n) H } \pm o(1) \nln
    \sk{c}{\in} & \EE[\mu]{\p{h(0) + h'(0)\frac{H}{n}}H} \pm c_1\EE[\mu]{\frac{H^2}{n^2}H } \pm c_2\pb(\overline{\calE}) \pm o(1) \nln
  {\in} & \EE[\mu]{\p{h(0) + h'(0)\frac{H}{n}}H} \pm c_1\EE[\mu]{\frac{H^2}{n^2}H } \pm o(1) \nln
   \sk{d}{\in} & \EE[\mu]{\p{h(0) + h'(0)\frac{H}{n}}H} \pm \mcal{O}(1/n) \pm o(1) \nln
   \sk{e}{\in} & \frac{d}{d\mu} \omega\p{\frac{d/n}{\mu}} \EE[\mu]{H^2}\pm o(1)  \nln
   =  & \frac{d}{d\mu} \omega\p{\frac{d/n}{\mu}} \mu\nth_a(1-\mu) \pm o(1)  \nln
    =  & \frac{d}{d\mu} \omega\p{\frac{d/n}{\mu}} (\eta'(\mu))^{-1}  \pm o(1)  
\label{eq:EuH}
 \end{align} 
   where $ c_2 = \max_{d,x \in \rp}\Omega(d,x),  \mbox{ and }  c_3 = \max_{x \in [-\delta, \delta]} h(x)$, and the $o(1)$ term does not depend on $d$.  Step $(a)$ is based on the Cauchy–Schwarz inequality, $(b)$ from the fact that $\pb\p{\overline{\calE}}$ converges to $0$ exponentially fast in $n$ by the Chernoff bound and that $\EE{H^2} = \mcal{O}(n)$, $(c)$ from the Taylor expansion in \eqref{eq:hfcttaylor}, and $(d)$ from the fact that $\abs{\EE{H^3}} = \mathcal{O}(n)$ as a result of $X$ being a Binomial random variable. Finally, step $(e)$ follows from the definition of $h$ in \eqref{eq:hfctdef}. 

Recall from \eqref{eq:Dconv0} that, conditional on $A=a$, $D/n$ concentrates on $d_a$ as $n\to \infty$. \eqref{eq:EuH} thus implies that 
\begin{align}
\lim_{n\to \infty}   \EE[\mu]{\omega\p{{D}/{X}}H }   = & \lim_{n\to \infty} \sum_{d\in \zp}  \EE[\mu]{\omega\p{{d}/{X}}H }   \PP{D=d |A=a}\nln
= & \frac{d}{d\mu} \omega\p{{d_a}/{\mu}} (\eta'(\mu))^{-1}
\label{eq:EuHrand}
\end{align}
Substituting \eqref{eq:ElDX} and \eqref{eq:EuHrand} into \eqref{eq:dmuQdX}, we obtain that 
\begin{equation}
\lim_{n \to \infty} \frac{d}{d\mu} \EE[\mu]{ \Omega \p{D, X}|A=a} = \frac{d}{d\mu} \omega\p{{d_a}/{\mu}}. 
\end{equation}
This proves Proposition  \ref{prop:qprimeconvg}. \qed

\subsection{Proof of Proposition \ref{prop:perturb_df}}
\label{app:prop:perturb_df}

Fix $a\in \calA$ and $n\in \N$.  
Denote by $\pi_{p, \zeta}$ the $\zeta$-perturbed payment distribution centered at $p$ (\ref{eq:zEta}).  We first prove that $ \cb{\frac{\partial }{\partial \zeta} \mu\nth_a(p, \zeta)}_{\zeta=0}=0. $  By \eqref{eq:balanceq}, $\mu\nth_a(p,\zeta)$ satisfies
\begin{equation}
 \frac{\partial^k}{\partial^k \zeta} \mu\nth_a(p,\zeta) = \frac{\partial^k}{\partial^k \zeta} \psi\nth_a(\mu\nth_a(p,\zeta), \pi_{p,\zeta}), \quad k \in \N.
 \label{eq:balanceEqkth}
 \end{equation}
It therefore suffices to evaluate the right-hand side of the above equation. To this end: 
\begin{align}
& \psi\nth_a(\mu\nth_a(p,\zeta), \pi_{p,\zeta})- \psi\nth_a(\mu\nth_a(p), \delta_p)  \nln
\sk{a}{=} & \frac{1}{2}\lt( \EE { f_{B_1} \lt( (p+\zeta) \,q\nth_a(\mu\nth_a(p, \zeta)) \rt) \bbar A=a }+ \EE{ f_{B_1} \lt( (p-\zeta) \,q\nth_a(\mu\nth_a(p, \zeta)) \rt)  \bbar A=a } \rt)  \nln
& - \EE{ f_{B_1} \lt( p \,q\nth_a(\mu\nth_a(p)))  \rt) \bbar A=a } \nln
= & \frac{1}{2}\lt( \EE{ f_{B_1} \lt( (p+\zeta) \,q\nth_a(\mu\nth_a(p, \zeta))  \rt)  \bbar A=a } - \EE{f_{B_1} \lt( p \,q\nth_a(\mu\nth_a(p)))  \rt) \bbar A=a } \rt) \nln
& +  \frac{1}{2}\lt( \EE{ f_{B_1} \lt( (p-\zeta) \,q\nth_a(\mu\nth_a(p, \zeta)) \rt) \bbar A=a}- \EE{ f_{B_1} \lt( p \,q\nth_a(\mu\nth_a(p)))  \rt)  \bbar A=a } \rt), 
\label{eq:diffPsi}
\end{align}
where $(a)$ follows from the definition of $\zeta$-perturbation (\eqref{eq:zEta})  and the independence of perturbations $\{\varepsilon_i\}_{i \in \N}$ from the rest of the system. 
Since both $f_{B_1}(\cdot)$ and $q\nth_a(\cdot)$ are bounded,  for the first term on the right-hand side of \eqref{eq:diffPsi}, it is not difficult to show using the dominated convergence theorem that   there exists $c>0$ such that\footnote{Notation: $x \in y \pm z \leftrightarrow x \in [y-z, y+z]$.}
\begin{align}
&  \EE{ f_{B_1} \lt( (p+\zeta) \,q\nth_a(\mu\nth_a(p, \zeta)) \rt) \bbar A=a} - \EE{ f_{B_1} \lt( p \,q\nth_a(\mu\nth_a(p)))  \rt)  \bbar A=a} \nln
\in &  \upsilon \zeta  \, \pm \, c\zeta^2, 
 \label{eq:EsigbDiff}
\end{align}
for all sufficiently small $\zeta$, where   $ \upsilon \bydef \cb{ \frac{\partial }{ \partial \zeta}  \EE{ f_{B_1} \lt( (p+\zeta) \,q\nth_a(\mu\nth_a(p,\zeta)) \rt) \bbar A=a}}_{\zeta=0}. $ Applying the same argument to the second term in \eqref{eq:diffPsi}, we have that there exists $c$, such that for all sufficiently small $\zeta$
\begin{equation}
\frac{\psi\nth_a(\mu\nth_a(p,\zeta), \zeta)- \psi\nth_a(\mu\nth_a(p,0), 0)}{\zeta}  \in \pm  \frac { c\zeta^2}{\zeta} = \pm c\zeta, 
\end{equation}
which further implies that
\begin{equation}
\cb{\frac{\partial }{\partial \zeta} \mu\nth_a(p, \zeta)}_{\zeta=0}=0.
\label{eq:dmudzeta}
\end{equation}

For the derivative of $q\nth_a(\mu\nth_a(p,\cdot))$,  note that by chain rule, we have
\begin{align}
\frac{\partial }{\partial \zeta} q\nth_a(\mu\nth_a(p,\zeta))  = (q\nth_a)'(\mu\nth_a(p,\zeta))\frac{\partial}{\partial \zeta}\mu\nth_a(p,\zeta). 
\label{eq:omegader1}
\end{align}
Since $\cb{\frac{\partial}{\partial \zeta}\mu\nth_a(p,\zeta)}_{\zeta=0} = 0$  by \eqref{eq:dmudzeta}, and  $ (q\nth_a)'(\mu\nth_a(p,\zeta))$ is finite  by Proposition \ref{prop:qprimeconvg},  we have that $\cb{\frac{\partial }{\partial \zeta}q\nth_a(\mu\nth_a(p, \zeta)))}_{\zeta=0}=0$. This proves the first claim of Proposition \ref{prop:perturb_df}.

For the second claim, define 
\begin{equation}
g_{\varepsilon_1}(\zeta) \bydef \p{p+\varepsilon_1\zeta}q\nth_a \p{\mu\nth_a(p,\zeta)}. 
\end{equation}
Applying the chain rule to \eqref{eq:balanceEqkth}, we have that
\begin{align}
\frac{\partial^2 }{\partial^2 \zeta}\mu\nth_a(p,\zeta) =&  \EE{\frac{\partial^2}{\partial^2 \zeta}f_{B_1}\lt(g_{\varepsilon_1}(\zeta)\rt) \bbar A=a} \nln
=&   \EE{f''_{B_1}(g_{\varepsilon_1}(\zeta))g'_{\varepsilon_1}(\zeta)^2   + f'_{B_1}(g_{\varepsilon_1}(\zeta))g''_{\varepsilon_1}(\zeta) \bbar A=a}
\label{eq:part2mu}
\end{align}
Note that 
\begin{align}
g'_{\varepsilon_1}(\zeta) = & \frac{\partial}{\partial \zeta} \left[  (p+\varepsilon \zeta)q\nth_a(\mu\nth_a(p,\zeta)) \right] \nln
= & p\frac{\partial}{\partial \zeta}q\nth_a(\mu\nth_a(p,\zeta)) +  \varepsilon_1 q\nth_a(\mu\nth_a(p,\zeta)) + \varepsilon_1\zeta\frac{\partial}{\partial\zeta}q\nth_a(\mu\nth_a(p,\zeta)),
\label{eq:g1par}
\end{align}
and
\begin{align}
g''_{\varepsilon_1}(\zeta) =  p\frac{\partial^2}{\partial^2 \zeta}q\nth_a(\mu\nth_a(p,\zeta)) + 2\varepsilon_1 \frac{\partial}{\partial \zeta}q\nth_a(\mu\nth_a(p,\zeta)) + \varepsilon_1\zeta \frac{\partial^2}{\partial^2 \zeta} q\nth_a(\mu\nth_a(p,\zeta)). 
\label{eq:g2par}
\end{align}
By chain rule, we have 
\begin{align}
& \cb{\frac{\partial^2 }{\partial^2 \zeta}q\nth_a(\mu\nth_a(p,\zeta)))}_{\zeta=0} \nln
= & (q\nth_a)''(\mu\nth_a(p,0)) \cb{\frac{\partial}{\partial \zeta}\mu\nth_a(p,\zeta)}_{\zeta=0}^2 +(q\nth_a)'(\mu\nth_a(p,0)) \cb{\frac{\partial^2}{\partial^2 \zeta}\mu\nth_a(p,\zeta)}_{\zeta=0}\nln
= & (q\nth_a)'(\mu\nth_a(p)) \cb{\frac{\partial^2}{\partial^2 \zeta}\mu\nth_a(p,\zeta)}_{\zeta=0}, 
\label{eq:omegader2}
\end{align}
where the last step follows from the fact that $\cb{\frac{\partial}{\partial \zeta}\mu\nth_a(p,\zeta)}_{\zeta=0}=0$.  Applying  \eqref{eq:omegader1} and \eqref{eq:omegader2} to \eqref{eq:g1par} and \eqref{eq:g2par}, we have 
\begin{equation}
g'_{\varepsilon_1}(0) =p\cb{\frac{\partial}{\partial \zeta}q\nth_a(\mu\nth_a(p,\zeta))}_{\zeta=0} +  \varepsilon_1 q\nth_a(\mu\nth_a(p,0)) + 0 = \varepsilon_1 q\nth_a(\mu\nth_a(p)), 
\end{equation}
and 
\begin{align}
g''_{\varepsilon_1}(0) =&  p\cb{\frac{\partial^2}{\partial ^2 \zeta}q\nth_a(\mu\nth_a(p,\zeta))}_{\zeta=0}+ 2\varepsilon_1 \cb{\frac{\partial}{\partial \zeta}q\nth_a(\mu\nth_a(p,\zeta))}_{\zeta=0}+ 0  \nln
=&  p (q\nth_a)'(\mu\nth_a(p)) \cb{\frac{\partial^2}{\partial^2 \zeta}\mu\nth_a(p,\zeta)}_{\zeta=0}.  
\end{align}
Substituting the expressions for $g'_{\varepsilon_1}(0)$ and $g''_{\varepsilon_1}(0)$ into \eqref{eq:part2mu}, we obtain:
\begin{align}
& \cb{\frac{\partial^2 }{\partial^2 \zeta}\mu\nth_a(p,\zeta)}_{\zeta=0} \nln
=&  \EE{f''_{B_1}(g_{\varepsilon_1}(0))\varepsilon_1^2q\nth_a(\mu\nth_a(p))^2 + f'_{B_1}(g_{\varepsilon_1}(0))p(q\nth_a)'(\mu\nth_a(p)) \cb{\frac{\partial^2}{\partial^2 \zeta}\mu\nth_a(p,\zeta)}_{\zeta=0}   \cond A=a  }\nln
= &  \EE{f''_{B_1}(g_{\varepsilon_1}(0))  \cond A=a}q\nth_a(\mu\nth_a(p))^2 
\nln & + \EE{f'_{B_1}(g_{\varepsilon_1}(0))  \cond A=a} p(q\nth_a)'(\mu\nth_a(p))\cb{\frac{\partial^2}{\partial^2 \zeta}\mu\nth_a(p,\zeta)}_{\zeta=0}, 
\end{align} 
where the last step follows from the fact that $\varepsilon_1 \in \{-1, 1\}$ and hence $\varepsilon_1^2=1$. After re-arrangement, the above equation yields
\begin{equation}
\cb{\frac{\partial^2 }{\partial^2 \zeta}\mu\nth_a(p,\zeta)}_{\zeta=0}  = \frac{ \EE{f''_{B_1}(g_{\varepsilon_1}(0))  \cond A=a }q\nth_a(\mu\nth_a(p))^2 }{1- \EE{f'_{B_1}(g_{\varepsilon_1}(0))  \cond A=a}(q\nth_a)'(\mu\nth_a(p))p }, 
\label{eq:muDerv2}
\end{equation}
and by \eqref{eq:omegader2}, we have
\begin{align}
& \cb{\frac{\partial^2 }{\partial^2 \zeta}q\nth_a(\mu\nth_a(p,\zeta),\zeta))}_{\zeta=0}  \nln
= & (q\nth_a)'(\mu\nth_a(p,0)) \cb{\frac{\partial^2}{\partial^2 \zeta}\mu\nth_a(p,0)}_{\zeta=0} \nln
=&  (q\nth_a)'(\mu\nth_a(p))\lt( \frac{ \EE{f''_{B_1}(g_{\varepsilon_1}(0)) \cond A=a}q\nth_a(\mu\nth_a(p))^2 }{1- \EE{f'_{B_1}(g_{\varepsilon_1}(0))  \cond A=a}(q\nth_a)'(\mu\nth_a(p))p} \rt). 
\label{eq:omegaDerv2}
\end{align}
Finally, we check the uniform boundedness of the second derivatives with respect to all $n$ and all sufficiently small $\zeta$. To show that $\cb{\frac{\partial^2 }{\partial^2 \zeta}\mu\nth_a(p,\zeta)}_{\zeta=0} $, note that $f'_{B_1}(\cdot)$ is non-negative and $(q\nth_a)'(\cdot)$ non-positive (Proposition \ref{prop:qprimeconvg}). Therefore, the term  $\EE{f'_{B_1}(g_{\varepsilon_1}(0))  \cond A=a}(q\nth_a)'(\mu\nth_a(p))p$ is non-positive. By \eqref{eq:muDerv2}, this implies the uniform boundedness of $\frac{\partial^2 }{\partial^2 \zeta}\mu\nth_a(p,\zeta)$. Note that by Proposition \ref{prop:qprimeconvg},  $(q\nth_a)'(\mu\nth_a(p))$ is non-positive and  bounded, and with \eqref{eq:omegaDerv2} this shows that  $\frac{\partial^2 }{\partial^2 \zeta}q\nth_a(\mu\nth_a(p,\zeta))$ is bounded for all $n$ and all sufficiently small $\zeta$. This proves the second claim and thus completes the proof of  Proposition \ref{prop:perturb_df}.

\subsection{Proof of Proposition \ref{prop:utilityderiv}}
\label{app:prop:utilityderiv}

Fix $a \in \calA$ and $n\in \N$.  Consider the case where the payment distributions where all potential suppliers are offered a fixed payment, $p$, i.e., $\pi = \delta_p$.  Recall from \eqref{eq:reward} and \eqref{eq:chia} that 
\begin{align}
\frac{d}{dp}u_a\nth(p) =& \frac{d}{dp}\frac{1}{n}\EE{R(D, \, T) - \sum_{i = 1}^n P_i Z_i S_i \Bbar A=a} \nln
= & \frac{d}{dp}\EE{\frac{1}{n}R(D, \, T) \Bbar A=a} -   \frac{d}{dp} \p{ \frac{1}{n}p \EE{ \sum_{i=1}^n Z_i  S_i} \Bbar A=a} \nln
= & \frac{d}{dp}\EE{\frac{1}{n}R(D, \, T)  \Bbar A=a} -   \frac{d}{dp} \p{p \EE{\frac{1}{n} \Omega(D,\,T)T }\Bbar A=a}, 
\label{eq:dchi}
\end{align}
where  $T \sim\text{Binomial}(\mu\nth_a(p), \, n)$. We have by the chain rule:
\begin{align}
\frac{d}{dp} \EE{R(D, \, T)}=& \EE{(\mu\nth_a)'(p)\cb{ \frac{d}{d\mu}\EE[\mu]{R(D, \, X) \bbar A=a } }_{\mu = \mu\nth_a(p)}}
\label{eq:dRdp}
\end{align}
where $X \sim \text{Binomial}(\mu, n)$, and similarly
\begin{align}
\frac{d}{dp} \EE{\Omega(D, \, T)T} =\EE{(\mu\nth_a)'(p)\cb{ \frac{d}{d\mu}\EE[\mu]{\Omega(D, \, X)T \bbar A=a} }_{\mu = \mu\nth_a(p)}}. 
\label{eq:dQTdp}
 \end{align} 
Using arguments essentially identical to that of Proposition \ref{prop:qprimeconvg}, we can show that for all $a\in \calA$ and $\mu>0$
\begin{align}
\label{eq:dRdp_Stein}
&\limn \frac{d}{d\mu}\EE[\mu]{\frac{1}{n}R(D, \, X ) \bbar A=a}  =   \frac{d}{d\mu} \p{r(d_a/\mu) \mu} \\
&\limn \frac{d}{d\mu}\EE[\mu]{\frac{1}{n}\Omega(D, \, X )X \bbar A=a} =   \frac{d}{d\mu} \p{\omega(d_a/\mu) \mu}, 
\label{eq:dQTdp_Stein}
\end{align}
where the limiting functions $\omega$ and $r$ are defined in Assumptions \ref{ass:Rtou}  and \ref{ass:Omega_tou}, respectively. 
Substituting \eqref{eq:dRdp_Stein} and \eqref{eq:dQTdp_Stein} into \eqref{eq:dRdp} and \eqref{eq:dQTdp}, respectively, and observing that 
\begin{equation}
\limn \EE{\Omega(D, T)T/n \bbar A=a} = \omega(d_a/\mu_a(p))\mu_a(p), 
\end{equation}
we have 
\begin{align*}
 &\limn \frac{d}{dp} \EE{\frac{1}{n}R(D, \, T) \Bbar A=a}  =  \mu_a'(p)   \cb{\frac{d}{d\mu}\p{r(d_a/\mu)\mu}}_{\mu = \mu_a(p)},  \\
 &\limn \frac{d}{dp}\p{p \EE{\frac{1}{n}\Omega(D, \, T)T   \Bbar A=a }} =   p \mu_a'(p) \cb{\frac{d}{d\mu}\p{\omega(d_a/\mu)\mu }}_{\mu=\mu_a(p)}   +\omega(d_a/\mu_a(p))\mu_a(p) ,  
\end{align*}
where $\mu_a(p) = \limn \mu\nth_a(p)$ is defined in Lemma \ref{lem:limit_mup}. 
\begin{align*}
 &  \limn \frac{d}{d p}  u_a\nth (p)  \nln
 =&   {\mu_a'(p)   \cb{\frac{d}{d\mu}\p{r(d_a/\mu)\mu}}_{\mu = \mu_a(p)}}   -  {p \mu_a'(p) \cb{\frac{d}{d\mu}\p{\omega(d_a/\mu)\mu }}_{\mu=\mu_a(p)}   +\omega(d_a/\mu_a(p))\mu_a(p) } \nln
=&   \mu_a'(p) \p{ \cb{\frac{d}{d\mu}\p{r(d_a/\mu)\mu}}_{\mu = \mu_a(p)} - p\cb{\frac{d}{d\mu}\p{\omega(d_a/\mu)\mu}}_{\mu = \mu_a(p)}} \nln
 & \ \ - \omega(d_a/\mu_a(p))\mu_a(p) \nln
 = &  \mu_a'(p)\left[ r\p{\frac{d_a}{\mu_a(p)}} - p \omega\p{\frac{d_a}{\mu_a(p)}} - \p{ r'\p{\frac{d_a}{\mu_a(p)}} - p \omega'\p{\frac{d_a}{\mu_a(p)}}}\frac{d_a}{\mu_a(p)}\right] \nln
& \ \   - \omega\p{\frac{d_a}{\mu_a(p)}}\mu_a(p) \nln
 = & u_a'(p), 
\end{align*}
where $u_a(\cdot)$ is defined in \eqref{eq:uconvg}. This recovers the desired result.
\qed

\subsection{Proof of Proposition \ref{prop:demincreasp}}
\label{app:prop:demincreasp}
By the chain rule, and the fact that $q_a(\mu)= \omega(d_a/\mu)$, we have that 
\begin{align}
\frac{d}{dp}\p{q_a(\mu_a(p))\mu_a(p)} = & q'(\mu_a(p))\mu_a'(p) \mu_a(p) + q_a(\mu_a(p)) \mu_a'(p) \nln
= & (\omega(d_a/\mu_a(p)) - \omega'(d_a/\mu_a(p))d_a/\mu_a(p) )\mu_a'(p)
\end{align}
Using the expression for $\mu_a'(p)$ (cf.~\eqref{eq:mup2}), it is not difficult to show that, as a result of the strong concavity of $f_a(\cdot)$ in the interval $(\underline{x}, \overline{x})$, we have that $\inf_{p \in (c_0, p)} \mu_a'(p) >0$.   Furthermore, using the same argument as in the proof of Proposition \ref{prop:weakconcave} and the fact that $\omega(\cdot)$ is strongly concave with $\omega(0)\geq 0$, we have that $\inf_{p\in (c_0, \gamma)} (\omega(d_a/\mu_a(p)) - \omega'(d_a/\mu_a(p))d_a/\mu_a(p) )> 0$. Together, this implies that $\inf_{p \in (c_0, \gamma)}\frac{d}{dp}\p{q_a(\mu_a(p))\mu_a(p)} > 0$, thus proving our claim. 
\qed

\subsection{Proof of Proposition \ref{prop:Cmu'1neg}}
\label{app:prop:Cmu'1neg}
Note that 
\begin{align}
(\mu_a'(p), 1)   \bfC (\mu_a'(p), 1) ^\intercal  = \omega(d_a/\mu_a(p))\frac{2d_a \mu_a'(p)}{\mu_a(p)^2}\p{p\mu_a'(p) - \mu_a(p)}. 
\end{align}
It therefore suffices to show that 
\begin{equation}
p\mu_a'(p) - \mu_a(p) \leq 0. 
\label{eq:mu'weakconv}
\end{equation}
From Lemma \ref{lem:limit_mup}, we have that 
\begin{align}
 \mu_a'(p) = \frac{\Delta_a(p) }{1-p\Delta_a(p) q_a'(\mu_a(p)) /q_a(\mu_a(p))}, 
 \label{eq:muprimea2}
\end{align}
where 
\begin{equation}
\Delta_a(p) = q_a(\mu_a(p))f'_a(pq_a(\mu_a(p)), 
\end{equation}
and 
\begin{equation}
\mu_a(p) = f_a(pq_a(\mu_a(p))). 
\end{equation}
Multiplying the left-hand side of \eqref{eq:muprimea2} by $p/\mu_a(p)$, we obtain 
\begin{align}
\frac{\mu_a'(p)p}{\mu_a(p)}= & \frac{p\Delta_a(p)/\mu_a(p)} {1-p\Delta_a(p) q_a'(\mu_a(p)) /q_a(\mu_a(p)) } \nln
= & \frac{f'_a(pq_a(\mu_a(p)))\p{pq_a(\mu_a(p))}}{f\p{pq_a(\mu_a(p))}} \cdot \frac{1}{1-p\Delta_a(p) q_a'(\mu_a(p)) /q_a(\mu_a(p)) } \nln
\leq & \frac{\tilde{f}'(pq_a(\mu_a(p)))\p{pq_a(\mu_a(p))}}{\tilde{f}\p{pq_a(\mu_a(p))}} \cdot \frac{1}{1-p\Delta_a(p) q_a'(\mu_a(p)) /q_a(\mu_a(p)) } \nln
\sk{a}{\leq}  & \frac{1}{1-p\Delta_a(p) q_a'(\mu_a(p)) /q_a(\mu_a(p)) } \nln
\sk{b}{\leq}  & 1, \nnb
\end{align}
where $(a)$ follows from Proposition \ref{prop:weakconcave} combined with the non-negativity and concavity of $\tilde{f}(\cdot)$, and $(b)$ from the fact that $q_a'(\mu) = - \omega(d_a/\mu)d_a/\mu^2\leq 0$. This proves \eqref{eq:mu'weakconv} and hence the proposition. 
\qed

\ifx \useplain\undefined
\end{APPENDICES}
\else
\end{appendix}
\fi

\end{document}

%% file: title_info/informs_title.tex
\TITLE{Experimenting in Equilibrium}
\RUNTITLE{Experimenting in Equilibrium}

\RUNAUTHOR{Wager and Xu}

\ARTICLEAUTHORS{%

\AUTHOR{Stefan Wager}
\AFF{Graduate School of Business\\ Stanford University\\
 \EMAIL{swager@stanford.edu}}  

\AUTHOR{Kuang Xu}
\AFF{Graduate School of Business\\ Stanford University\\
 \EMAIL{kuangxu@stanford.edu}}  
 
} 

%% file: title_info/plain_title.tex
\usepackage{datetime}
\newdateformat{monthyeardate}{%
 \monthname[\THEMONTH], \THEYEAR}

\author{Stefan Wager \\
Graduate School of Business\\
Stanford University\\
\texttt{swager@stanford.edu} 
\and
Kuang Xu \\
Graduate School of Business\\
Stanford University\\
 \texttt{kuangxu@stanford.edu}}
\date{Draft Version \monthyeardate\today}
\title{Experimenting in Equilibrium}